\pdfoutput=1
\RequirePackage{ifpdf}
\ifpdf % We are running pdfTeX in pdf mode
\documentclass[pdftex]{sigma}
\else
\documentclass{sigma}
\fi

\numberwithin{equation}{section}

\newtheorem{Theorem}{Theorem}[section]
\newtheorem{Corollary}[Theorem]{Corollary}
\newtheorem{Lemma}[Theorem]{Lemma}
\newtheorem{Proposition}[Theorem]{Proposition}
 { \theoremstyle{definition}
\newtheorem{Definition}[Theorem]{Definition}

\newtheorem{Example}[Theorem]{Example}
\newtheorem{Remark}[Theorem]{Remark}
\newtheorem{Notation}[Theorem]{Notation}}

\newcommand{\defeq}{\overset{\scriptstyle \mathrm{def.}}{=}}
\newcommand{\C}{{\mathbb C}}
\newcommand{\R}{{\mathbb R}}
\newcommand{\Z}{{\mathbb Z}}
\newcommand{\Q}{{\mathbb Q}}

\newcommand{\M}{{\mathcal{M}}}

\newcommand{\w}{\mathbf{w}}

\newcommand{\vb}{\mathbf{v}}

\newcommand{\rank}{\mathop{\mathrm{rank}}}
\newcommand{\Hom}{\mathop{\mathrm{Hom}}}
\newcommand{\End}{\mathop{\mathrm{End}}}

\newcommand{\CP}{\operatorname{\C {\mathbb{P}}}}

\newcommand{\GL}{\operatorname{\mathrm GL}}
\newcommand{\SL}{\operatorname{\mathrm SL}}

\newcommand{\Hilb}[2]{{#1}^{\lbrack{#2}\rbrack}} % Hilbert scheme of points
\newcommand{\Hilbn}[1]{\Hilb{#1}{n}}

\newcommand{\shfO}{\mathcal O} % sheaf
\newcommand{\Supp}{\operatorname{Supp}} % support
\newcommand{\mapright}[1]{%
 \smash{\mathop{%
 \hbox to 1cm{\rightarrowfill}}\limits^{#1}}}
\newcommand{\mapleft}[1]{%
 \smash{\mathop{%
 \hbox to 1cm{\lefttarrowfill}}\limits^{#1}}}

\def\Young#1{\vbox{\smallskip\offinterlineskip
 \halign{&\vbox{##}\kern-\Thickness\cr #1}}}

\newdimen\Squaresize \Squaresize=12pt
\newdimen\Thickness \Thickness=.5pt
\newdimen\Correction \Correction=7pt

\def\Gauche#1{\hbox{\vrule width \Thickness
 \vbox to \Squaresize{\vss
 \hbox to \Squaresize{\hss#1\hss}
 \vss}
 \unskip\kern\Thickness}
 \kern-\Thickness}

\def\Vide#1{\hbox{
 \vbox to \Squaresize{\vss
 \hbox to \Squaresize{\hss#1 \hss}\vss}
 \hskip-\Correction}
 \kern-\Thickness}

\def\Droite#1{\hbox{\kern\Thickness
 \vbox to \Squaresize{\vss
 \hbox to \Squaresize{\hss#1\hss}
 \vss}
 \unskip\vrule width \Thickness}
 \kern-\Thickness}

\def\Haut#1{\hbox{\kern-\Thickness
 \vbox to \Squaresize{\hrule height \Thickness\vss
 \hbox to \MyDim{\hss#1\hss}
 \vss}
 \unskip}
 \kern-\Thickness}

\def\Bas#1{ \hbox{\kern-\Thickness
 \vbox to \Squaresize{\vss
 \hbox to \MyDim{\hss#1\hss}

 \vss\hrule height\Thickness}
 \unskip}
 \kern-\Thickness}

\def\Carre#1{\hbox{\vrule width \Thickness
 \vbox to \Squaresize{\hrule height \Thickness\vss
 \hbox to \Squaresize{\hss#1\hss}
 \vss\hrule height\Thickness}
 \unskip\vrule width \Thickness}
 \kern-\Thickness}

\def\Ca{\Carre{}}

\def\Box#1{\Carre{$#1$}}

\begin{document}

\allowdisplaybreaks

\newcommand{\arXivNumber}{math.AG/0510455}

\renewcommand{\PaperNumber}{052}

\FirstPageHeading

\ShortArticleName{A Combinatorial Study on Quiver Varieties}

\ArticleName{A Combinatorial Study on Quiver Varieties}

\Author{Shigeyuki FUJII~$^\dag$ and Satoshi MINABE~$^\ddag$}

\AuthorNameForHeading{S.~Fujii and S.~Minabe}

\Address{$^\dag$~Accenture Strategy, 107-8672 Tokyo, Japan}
\EmailD{\href{mailto:shigeyuki.fujii@accenture.com}{shigeyuki.fujii@accenture.com}}

\Address{$^\ddag$~Department of Mathematics, Tokyo Denki University, 120-8551 Tokyo, Japan}
\EmailD{\href{mailto:minabe@mail.dendai.ac.jp}{minabe@mail.dendai.ac.jp}}

\ArticleDates{Received January 13, 2017, in f\/inal form June 30, 2017; Published online July 06, 2017}

\Abstract{This is an expository paper which has two parts. In the f\/irst part, we study quiver varieties of af\/f\/ine $A$-type from a combinatorial point of view. We present a combinatorial method for obtaining a closed formula for the generating function of Poincar\'e polynomials of quiver varieties in rank~1 cases. Our main tools are cores and quotients of Young diagrams. In the second part, we give a brief survey of instanton counting in physics, where quiver varieties appear as moduli spaces of instantons, focusing on its combinatorial aspects.}


\Classification{14C05; 14D21; 05A19; 05E10}

\section{Introduction}\label{sec:intro}
Quiver varieties, which were introduced by Nakajima \cite{Nakajima2}, stand at the corner of a beautiful interplay between combinatorics, geometry, representation theory, and mathematical physics. This expository paper has two purposes. The f\/irst one is to study quiver varieties of type~$A_{l-1}^{(1)}$ from a combinatorial viewpoint (Sections~\ref{sec:s3}, \ref{sec:s2} and~\ref{sec:enu}). The second is to give a brief survey of instanton counting in physics, especially results of Maeda et al.~\cite{Nakatsu3, Nakatsu1, Nakatsu2}, and to discuss combinatorial aspects of them (Section~\ref{section5}).

The goal in the f\/irst part is to compute Poincar\'e polynomials of quiver varieties in rank~$1$ case via a combinatorial method. Note that Poincar\'e polynomials of quiver varieties have already been obtained by Nakajima (see Remark~\ref{rem:nakajima}), and our method just conf\/irms his result. However, we hope that our combinatorial viewpoint gives a new insight into the subject. Our strategy is as follows. There is natural torus actions on quiver varieties. Then by localization principle, we can compute topologies of them by studying the f\/ixed point set. We see that they further reduce to enumeration of Young diagrams. More importantly, we compute the Poincar\'e polynomial of a quiver variety not one by one, but a \textit{generating function} of them. For that purpose, we use some combinatorial machineries, which are \textit{cores} and \textit{quotient} of Young diagrams. They arise from a division algorithm for Young diagrams analogous to that for integers and play important roles in our discussions. They also appear in studies of supersymmetric gauge theories by Maeda et al.~\cite{Nakatsu3, Nakatsu1, Nakatsu2} (in slightly dif\/ferent terminologies). The purpose of the second part is to explain it from a mathematical viewpoint. We hope that our exposition is also useful for physics oriented readers.

We introduce quiver varieties of type $A_{l-1}^{(1)}$ as follows. Let $\Gamma$ be a cyclic group of order~$l$. Then $\Gamma$ acts on the framed moduli space $\M_{r}(n)$ of torsion free sheaves $E$ on $\CP^2$ with rank $r$, $c_2 = n$, where the framing is a trivialization at the line at inf\/inity $l_{\infty}$. We consider its $\Gamma$-f\/ixed point components,
\begin{gather}
	\M_{r} (n)^{\Gamma} = \bigsqcup_{\vb \atop |\vb|=n} \M_{\w} (\vb),\label{eq:decomp}
\end{gather}
where $\w$ and $\vb$ are isomorphism classes of $\Gamma$-modules. These $\M_{\w}(\vb)$'s are \textit{quiver varieties of type} $A_{l-1}^{(1)}$. Note that this def\/inition, which is taken from \cite{Nakajima3}, is slightly dif\/ferent from the original one. The original def\/inition was based on the so-called ADHM description of instantons on so-called ALE spaces~\cite{KN}.

The space $\M_{r}(n)$ was also introduced by Nakajima~\cite{Nakajima8}, as a resolution of the moduli space of instantons on~$\R^4$. Then quiver varieties can be seen as resolutions of the moduli spaces of instantons on $\C^2/\Gamma$. It turns out that generating functions of Euler characteristics of quiver varieties can be identif\/ied with the partition functions of certain supersymmetric f\/ield theory, which is important in mathematical physics. This is one of our motivations to study quiver varieties. We summarize some basic facts of instanton counting in Section~\ref{section5}. Here we should mention a beautiful work by physicists Fucito~et~al.~\cite{FMP}. They studied instanton moduli spaces on ALE spaces by a combinatorial method. The present work grow out of an attempt at rigorous understanding of their idea mathematically. In this paper, we give a precise mathematical formulation of the idea in~\cite{FMP}, and complete proofs to the results announced there.

Let us explain more detail about our analysis in the f\/irst part. There is a natural action of $(r+2)$-dimensional torus $T$ on a quiver variety~$\M_{\w}(\vb)$, induced from that on $\M_{r}(n)$. We study their f\/ixed point sets. It can be shown that f\/ixed points on a quiver variety $\M_{\w}(\vb)$ are parametrized by a certain set $\mathcal{P}_{\w}(\vb)$ of Young diagrams, which is determined by the {\it coloring of Young diagrams} (see Section~\ref{sec:s3} for def\/inition). Then localization theorem tells us that we can compute Euler characteristics of quiver varieties via enumeration of $\mathcal{P}_{\w}(\vb)$. We can also compute Poincar\'e polynomials by studying tangent spaces at f\/ixed points. As mentioned above, what we want to compute is the generating function of them with respect to $\vb$:
\begin{gather*}
\mathcal{Z}_{\w} (\mathfrak{t},\mathfrak{q},\vec{\mathfrak{r}}) \defeq \sum_{m\geq 0 \atop \vb \in (\Z_{\geq 0})^{l}}\mathrm{b}_m ( \M_{\w}(\vb))\mathfrak{t}^m \mathfrak{q}^{|\vb|} \prod_{i=0}^{l-1} \mathfrak{r}_{i}^{v_i} .
\end{gather*}
The problem is that how to get a closed formula for $\mathcal{Z}_{\w} (\mathfrak{t},\mathfrak{q},\vec{\mathfrak{r}})$. In Section~\ref{sub:euler}, we give a~closed formula for $\mathcal{Z}_{\w} (\mathfrak{t},\mathfrak{q},\vec{\mathfrak{r}})$ when $\w$ is a class $\mathbf{e}_{j}$ of $1$-dimensional $\Gamma$-module. First, we show a remarkable factorization property of the generating function of the Poincar\'e polynomials in Theorem~\ref{thm:t2}. It decomposes into the product of the contributions of \textit{cores} and those of \textit{quotients} of Young diagrams in $\mathcal{P}_{\mathbf{e}_{j}}(\vb)$:
\begin{gather*}
	\mathcal{Z}_{\mathbf{e}_{j}} (\mathfrak{t},\mathfrak{q},\vec{\mathfrak{r}}) =\mathcal{Z}^{\rm{quot}} (\mathfrak{t},\mathfrak{q},\vec{\mathfrak{r}}) \mathcal{Z}_{\mathbf{e}_{j}}^{\rm{core}} (\mathfrak{q},\vec{\mathfrak{r}}).
	\end{gather*}
Next, we compute each part separately.
Computations of the core part $\mathcal{Z}_{\mathbf{e}_{j}}^{\rm{core}} (\mathfrak{q},\vec{\mathfrak{r}})$ is purely combinatorial. On the other hand, we use geometry of Hilbert scheme of points, to get the quotients part $\mathcal{Z}^{\rm{quot}} (\mathfrak{t},\mathfrak{q},\vec{\mathfrak{r}})$. Then we obtain a closed formula for $\mathcal{Z}_{\mathbf{e}_{j}} (\mathfrak{t},\mathfrak{q},\vec{\mathfrak{r}})$ in Theorem~\ref{thm:t1}, where $\mathcal{Z}_{\mathbf{e}_{j}} (\mathfrak{t},\mathfrak{q},\vec{\mathfrak{r}})$ is expressed as a ratio of the Riemann theta function and the Dedekind eta function. This result is due to Nakajima, as we mentioned earlier.

The paper is structured as follows. We devote Section~\ref{sec:s3} to combinatorial preliminaries. The most important notions in Section~\ref{sec:s3} are \textit{cores} and \textit{quotients} of Young diagrams. In Section~\ref{sec:s2}, we introduce quiver varieties and torus actions on them. We study their f\/ixed point sets. In Section~\ref{sec:enu}, we compute two kinds of global invariants on quiver varieties. In Section~\ref{sec:integral}, we consider `equivariant volumes' of quiver varieties, which are def\/ined in Section~\ref{sub:vol}. In Section~\ref{sub:euler}, we study generating functions of Poincar\'e polynomials and Euler characteristics of quiver varieties, by using combinatorics prepared in Section~\ref{sec:s3}. In Section~\ref{section5}, we summarize developments of instanton counting in physics (up to 2005) and discuss their combinatorial aspects. In Appendix~\ref{appendixA}, we collect some facts about Hilbert scheme of points which are used in this paper.

\subsection*{Note added in 2017}
The original version of this paper was submitted to the arXiv in 2005. In this version, we have updated the bibliography and added references \cite{Cirafici-Szabo,Gyenge2,Gyenge1,GNS,Szabo} in order to ref\/lect some recent developments in the subject. We brief\/ly comment on them.

(i) The generating function of Euler characteristics of quiver varieties is re-considered in \cite{Gyenge1}. The argument in loc.\ cit.\ is also combinatorial but dif\/ferent from ours. The method of loc.\ cit.\ is used in \cite{Gyenge2, GNS} to compute the generating functions of Euler characteristics of Hilbert schemes of points of the singular surface of type $A$ or $D$.

(ii) The study of instanton counting in gauge theories is still an active f\/ield of research in relation to geometric representation theory. For recent advances, see, e.g., \cite{Cirafici-Szabo, Szabo} and references therein. In \cite[Section~2]{Szabo}, results of this paper are nicely recast into the framework of instanton counting on quotient stacks $[\mathbb{C}^2/\mathbb{Z}_l]$. Their relations to instanton counting on ALE spaces (i.e., minimal resolutions of $\mathbb{C}^2/\mathbb{Z}_l$) are also clarif\/ied in loc.\ cit.

\section{Preliminaries on combinatorics}\label{sec:s3}
In this section, we f\/irst f\/ix some notations on Young diagrams which are used through the paper. We also introduce the notion of \textit{quotients} and \textit{cores} for Young diagrams. Our basic references in this section are \cite[Section~2.7]{GA} and \cite[Section~ 3]{Olsson} (see also \cite[Exercise~7.59]{Stanley} and \cite[I.1, Example~8]{Macdonald}).

\subsection{Notations on Young diagrams}

\subsubsection{Young diagrams}
A partition is a sequence of non-increasing non-negative integers $\lambda = (\lambda_1 ,\ldots , \lambda_m)$. The corresponding Young diagram is a collection of rows of square boxes which are left adjusted and with~$\lambda_i$ boxes in the $i$th row for $i=1, \ldots ,m$. We identify partitions with their corresponding Young diagrams.

\begin{Definition}\quad
\begin{enumerate}\itemsep=0pt
\item[(i)] The {\it{weight}} $|Y|$ of a Young diagram $Y$ is the number of boxes in $Y$. In terms of the parti\-tion~$\lambda$ corresponding to $Y$,	$|Y| = \sum_i \lambda_i$.
\item[(ii)] The {\it{total weight}} $|\underline{Y}|$ of an $l$-tuple of Young diagrams $\underline{Y}=(Y_{1}, \ldots, Y_{l} )$ is def\/ined by $|\underline{Y}| \defeq \sum\limits_{k=1, \ldots, l} |Y_{k}|$.
\end{enumerate}
\end{Definition}

\begin{Notation}\quad
\begin{enumerate}\itemsep=0pt
\item[(i)] Let $\mathcal{P}$ be the set of Young diagrams. Then we have
\begin{gather*}
\mathcal{P}=\bigsqcup_{n \geq 0} \mathcal{P}(n),
\end{gather*}
where $\mathcal{P}(n) = \{Y\in \mathcal{P} \,|\, |Y|=n \}$.
\item[(ii)] Let $\mathcal{P}_{l}$ be the set of $l$-tuples of Young diagrams. Then we have
\begin{gather*}
\mathcal{P}_{l} =\bigsqcup_{n \geq 0} \mathcal{P}_{l}(n),
\end{gather*}
where $\mathcal{P}_{l}(n) = \{\underline{Y}\in \mathcal{P}_{l} \,|\, |\underline{Y}|=n\}$.
\end{enumerate}
\end{Notation}

\subsubsection{Functions on Young diagrams}
For a Young diagram $Y$, we def\/ine $l_{Y} (s)$ and $a_{Y} (s)$ as
\begin{gather*}
l_{Y} (s) \defeq {\lambda}_h -k, \qquad a_{Y} (s) \defeq {\lambda}'_k -h,
\end{gather*}
where $s=(h,k) \in (\Z_{>0})^2$, ${\lambda}_h$ is the length of the $h$-th row of $Y$, and ${\lambda}'_k$ is the length of the $k$-th column of $Y$.
\begin{Definition}\label{def:App1}\quad
\begin{enumerate}\itemsep=0pt
\item[(i)] The box $s \in Y$ sitting on the $i$-th row and the $j$-th column is called $(i,j)$-component of $Y$ and denoted by $s=(i,j)$.
\item[(ii)] For $s \in Y$, we def\/ine the {\it{leg length}} by $l(s)\defeq l_Y (s)$, the {\it{arm length}} by $a(s)\defeq a_Y (s)$, and the {\it{hook length}} by $h(s) \defeq l(s)+a(s)+1$.
\end{enumerate}
\end{Definition}

In the later sections, we need the following generalization of the def\/inition of hook length.

\begin{Definition}
Take two Young diagrams $Y_{\alpha}$ and $Y_{\beta}$. For $s\in (\Z_{>0})^2$ we def\/ine
\begin{gather*}
h_{Y_{\alpha},Y_{\beta}} (s) \defeq l_{Y_{\beta}} (s) +a_{Y_{\alpha}} (s) +1.
\end{gather*}
We call $h_{Y_{\alpha},Y_{\beta}} (s)$ \textit{relative hook length of} $s$ \textit{with respect to} $Y_{\beta}$ when $s \in Y_{\alpha}$.
\end{Definition}

\begin{Example}
In Fig.~\ref{h}, the Young diagram drawn by dotted lines is $Y_\alpha$ and the Young diagram drawn by solid lines is~$Y_\beta$.
The number of diamonds (resp.\ stars) is $l_{Y_{\beta}} (s)$ (resp.~$a_{Y_{\alpha}} (s)$). Then we have $h_{Y_\alpha ,Y_\beta} (s)=7$.
\begin{figure}[h!]\centering
\includegraphics[width=2.5cm,keepaspectratio]{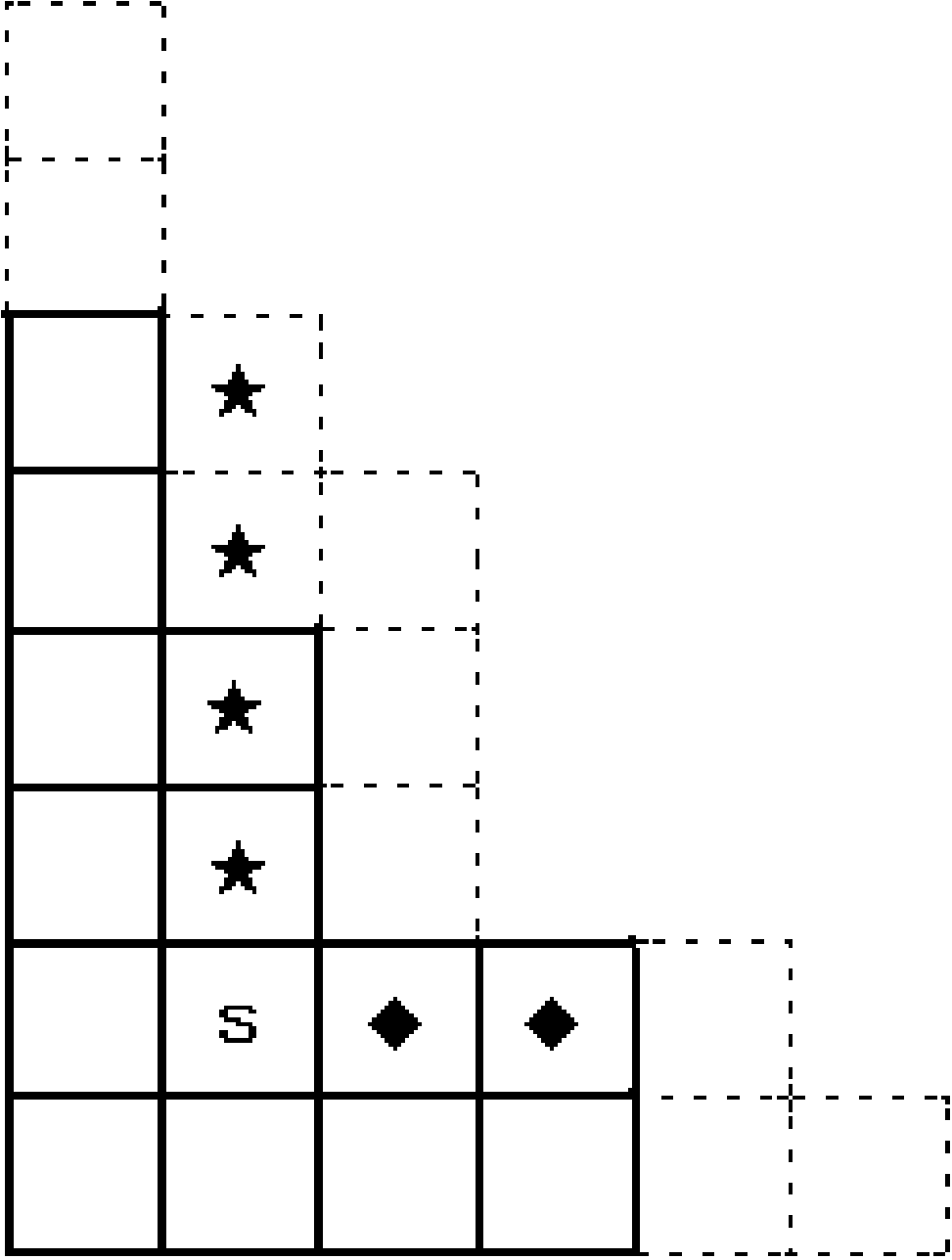}
\caption{$Y_\alpha$ (dotted) and $Y_\beta$ (solid).}\label{h}
\end{figure}
\end{Example}

\subsubsection{Removing a hook}
\begin{Definition}
For $s=(h,k) \in Y$, we def\/ine the {\it{hook}} of $s$ as
\begin{gather*}
H_s (Y) \defeq \bigl\{ (i,j)\in Y \,|\, i=h , \, j \geq k \ \mbox{or} \ j=k ,\ , i > h\bigr\}.
\end{gather*}
\end{Definition}
By def\/inition, $h(s) = \# H_s (Y)$. We call a hook with $n$ boxes a \textit{hook of length}~$n$.
\begin{Definition}
For $s=(h,k) \in Y$, we def\/ine the {\it{rim}} of $s$ as
\begin{gather*}
	R_s (Y) \defeq \bigl\{ (i,j)\in Y \,|\, (i+1,j+1)\notin Y \bigr\} \cap \bigl\{ (i,j)\in Y \,|\, i\geq h , \, j \geq k \bigr\}.
\end{gather*}
\end{Definition}
Then we have
\begin{gather*}
	\# R_s (Y) = \# H_s (Y) = h(s).
\end{gather*}

By def\/inition, removing a hook $H_s (Y)$ from $Y$ means that removing all boxes in~$R_s (Y)$ from~$Y$.
\begin{Example}
Let $Y$ be a Young diagram given by a partition $(5,4,3,2,1)$. We remove $H_{(2,2)} (Y)$ from $Y$ as follows:
\begin{gather*}\label{ex:rim_hook}
	\begin{array}{@{}c}
	\Young{\Box{}\cr
	 \Box{}&\Box{\star }\cr
	 \Box{}&\Box{\star}&\Box{\star}\cr
	 \Box{}&\Box{ }&\Box{\star }&\Box{\star }\cr	
 \Ca{}&\Box{}&\Box{}&\Box{}&\Box{}\cr}
 \end{array}
 \longmapsto
 \begin{array}{@{}c}
	\Young{\Box{}\cr
	 \Box{}\cr
	 \Box{}\cr
	 \Box{}&\Box{}\cr	
 \Ca{}&\Box{}&\Box{}&\Box{}&\Box{}\cr}
 \end{array}
\end{gather*}
\end{Example}

\subsubsection{Coloring of Young diagrams}
\begin{Definition}
Let $Y \in \mathcal{P}$. Let us f\/ix an integer $l \geq 2$. Assume that we are given a $\gamma \in \Z / l\Z$. Then the $l$-residue of a box $s \in Y$ is def\/ined by the following rule: if $s$ is the $(i,j)$-component of $Y$, then
\begin{gather*}
\operatorname{res}(s) = \gamma +i -j +l\Z \in \Z / l\Z.
\end{gather*}
We call this assignment of an element of $\Z / l\Z$ to each box in $Y$ as $\gamma$\textit{-coloring} of a Young diag\-ram~$Y$.
\end{Definition}
We extend this coloring to tuples of Young diagrams.
\begin{Definition}
Let $\underline{Y} = (Y_{1}, \ldots , Y_{r}) \in \mathcal{P}_{r}$. Assume that we are given $\gamma _1 , \ldots ,\gamma _r \in \Z /l\Z$. If a~box $s \in \underline{Y}$ is a $(i,j)$-component of $Y_{k}$, then we color $s$ with
\begin{gather*}
\operatorname{res}(s) = \gamma _k +i-j+l\Z \in \Z/l\Z.
\end{gather*}
We call this assignment as a $(\gamma_1 , \ldots , \gamma_r)$-\it{coloring} of an $r$-tuple of Young diagrams $\underline{Y}$.
\end{Definition}
\begin{Example}
Take $r=2$, $l=3$, and $(\gamma_1 , \gamma_2 )=(0,2) \in (\Z /3\Z)^2$. Let $\underline{Y} =(Y_1 , Y_2 ) \in \mathcal{P}_2$ which is given by a~pair of partitions $(\lambda_1 , \lambda_2 )=( (4,3,1) , (3,2,1,1) )$. Then $(\gamma_1 , \gamma_2 )$-coloring of~$\underline{Y}$ is given as follows:
\begin{gather*}
	\left(
	\begin{matrix}
 	\Young{\Box{1}\cr
	 \Box{2}&\Box{0}&\Box{1}\cr
	 \Box{0}&\Box{1}&\Box{2}&\Box{0}\cr} \, \, , \, \,
	\Young{\Box{2}\cr
	 \Box{0}\cr
 \Box{1}&\Box{2}\cr
	 \Box{2}&\Box{0}&\Box{1}\cr}
	 \end{matrix}
	 \right)
\end{gather*}
\end{Example}

\subsection{The Maya diagram}
\begin{Definition} A \textit{Maya diagram} is a sequence $\left\{ \mathbf{m} (\nu) \right\}_{\nu \in \Z }$ which consists of $0$ or $1$ and satisf\/ies the following property: there exist $N,M \in \Z$ such that for all $\nu >N$ (resp. $\nu <M$), $\mathbf{m} (\nu )=1$ (resp. $\mathbf{m} (\nu )=0$).
\end{Definition}
It is well-known that there is a one-to-one correspondence between the set of Maya diagrams and the set of Young diagrams. We explain how to associate a Young diagram to a Maya diagram. Let us place Young diagrams on the $(x,y)$-plane by the following way: the bottom-left corner of the Young diagram is at origin $(0,0)$ and a box in the Young diagram is an unit square. We call an upper-right borderline of $\{x\mbox{-axis}\} \cup \{ y\mbox{-axis} \} \cup Y$ the {\it{extended borderline}} of $Y$ and denote it by $\partial Y$. The line def\/ined by $y=x$ is called the {\it{medium}}.

\begin{figure}[h!]\centering
\includegraphics[width=4.5cm,keepaspectratio]{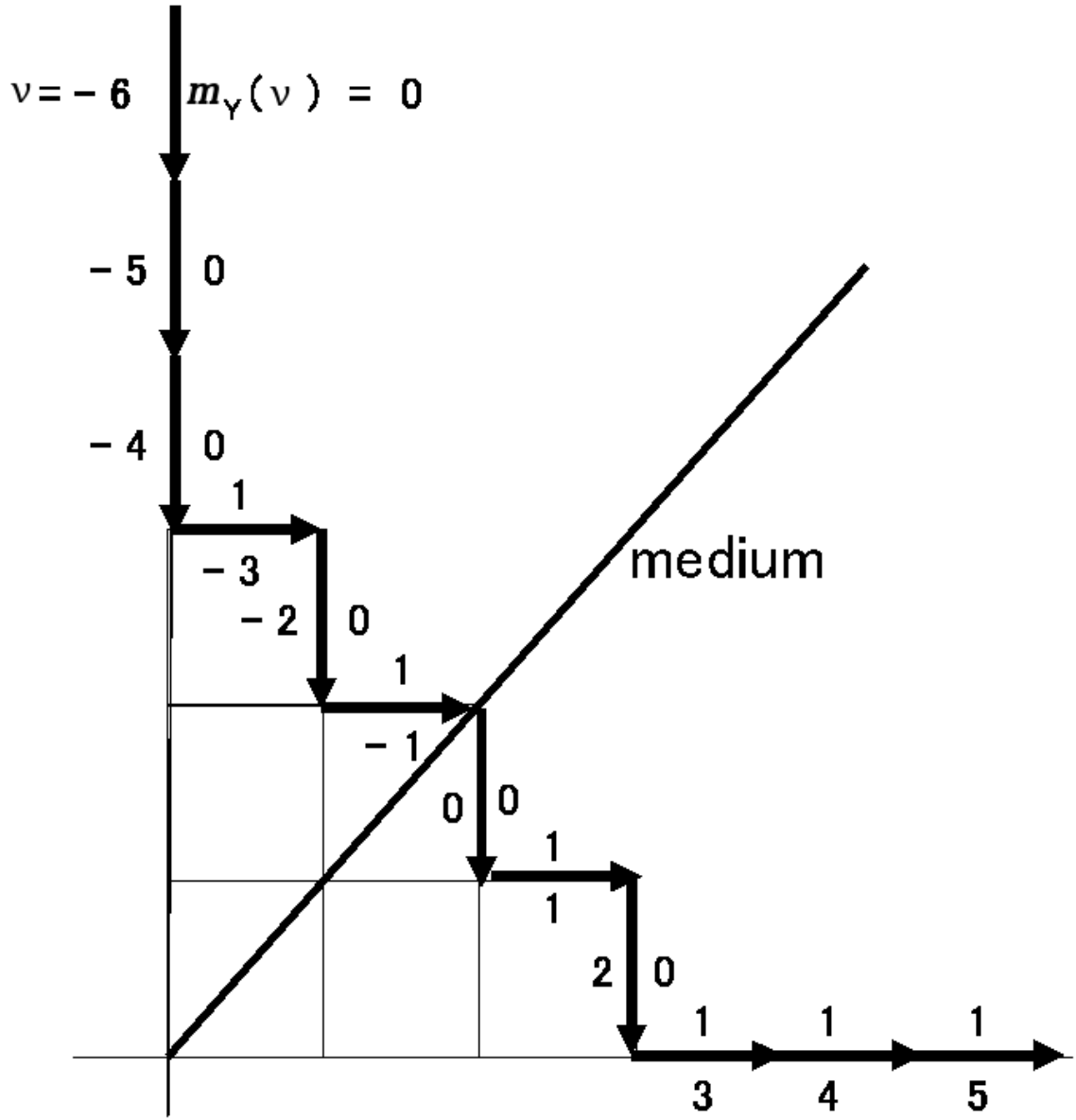}
\caption{Young diagram and Maya diagram.}\label{m}
\end{figure}

\begin{Example}
Let $Y$ be the Young diagram corresponding to the partition $(3,2,1)$. In Fig.~\ref{m}, $\partial Y$ is shown by a bold line.
\end{Example}

The Maya diagram $ \{ \mathbf{m}_Y (\nu) \}_{\nu \in \Z } $ corresponding to a Young diagram $Y$ is def\/ined as follows. We give the direction to the extended borderline $\partial Y$ of $Y$ as in Fig.~\ref{m}. Then each edge of the extended borderline is numbered by $\nu \in \Z$, if we set the edge which is located at the next to the medium to be~$0$. Next we encode a edge $\downarrow$ (resp.~$\rightarrow$) to $0$ (resp.~$1$). By this way, we have a $0/1$ sequence $\{\mathbf{m}_Y (\nu ) \}_{\nu \in \Z}$, where each $\mathbf{m}_Y (\nu )$ corresponds to a edge of $\partial Y$. Concretely, $ \{ \mathbf{m}_Y (\nu ) \}_{\nu\in \Z}$ is given by
\begin{gather*}
\mathbf{m}_Y (\nu )=
\begin{cases}
0 &\mathrm{if} \ \nu=\lambda_i-i,\\
1& \mathrm{else}.
\end{cases}
\end{gather*}
Or, equivalently, it is given by
\begin{gather*}
\mathbf{m}_Y (\nu )=
\begin{cases}
1 & \mathrm{if} \ \nu=i-1-\lambda_i',\\
0& \mathrm{else}.
\end{cases}
\end{gather*}
Here $\lambda_i$ (resp. $\lambda_i'$) is the length of the $i$-th row (resp.\ $i$-th column) of~$Y$.

\begin{Notation}\label{not:1}
For brevity, we denote an inf\/inite sequence of $0$'s (resp.~$1$'s) by $\underline{0}$ (resp.~$\underline{1}$).
\end{Notation}
\begin{Example}
The Maya diagram for the Young diagram in Fig.~\ref{m} is
\begin{gather*}
	\dots ,0,0,1,0,1 \,| \, 0,1,0,1,1, \dots.
\end{gather*}
We express this as $\underline{0} , 1,0,1 \,| \, 0,1,0, \underline{1} $. Here $|$ represents the position of the medium. This means that the medium's neighbor on the right is the $\mathbf{m}_Y (0)$.
\end{Example}

\begin{Remark}
Since the number of edges $\rightarrow$ before the medium is same as the number of edges $\downarrow$ after the medium, we have
\begin{gather*}
	\# \bigl\{ \mathbf{m}_Y (\nu )=1 \,|\,\nu <0\bigr\} =\# \bigl\{ \mathbf{m}_Y (\nu )=0 \,|\,\nu \geq 0\bigr\} .
\end{gather*}
\end{Remark}
\begin{Remark}[removing a hook in terms of the Maya diagrams] Fix $(x, y) \in Y$ whose hook length is $n$. Let $Y^{(x, y)}$ be the Young diagram obtained by
removing the hook of length $n$ corresponding to $(x,y)$ from $Y$. Then the Maya diagram $\{ \mathbf{m}_{Y^{(x, y)}} (\nu )\}_{\nu \in \Z }$ of $Y^{(x, y)}$ is given by
\begin{gather*}
	\mathbf{m}_{Y^{(x, y)}} (\nu )=
	\begin{cases}
	 \mathbf{m}_Y (\lambda_y-y) &\mbox{if} \ \nu = x-1-\lambda_x',\\
	 \mathbf{m}_Y (x-1-\lambda_x') &\mbox{if} \ \nu = \lambda_y-y, \\
		\mathbf{m}_Y (\nu ) &\mbox{else}.
	\end{cases}
\end{gather*}
\end{Remark}
\begin{Example}
Let $Y$ be the Young diagram shown in Fig.~\ref{m}. Let us remove the hook of length~$5$ corresponding to the box $(x,y)=(1,1)$ from~$Y$. Then the Maya diagram of $Y^{(1,1)}$ is obtained by exchanging $\mathbf{m}_Y (-3)=1$ and $\mathbf{m}_Y (2)=0$ as follows:
\begin{gather*}
\underline{0} ,\widehat{1} ,0,1 \,| \, 0,1,\widehat{0} ,\underline{1} \longmapsto \underline{0} ,\widehat{0} ,0,1 \,| \, 0,1,\widehat{1} ,\underline{1} ,
\end{gather*}
where $\widehat{}$ represents locations of exchanges.
\end{Example}

\subsection{Quotients and cores for Young diagrams}\label{subsec:2.4}

\subsubsection{Quotients}
For each $0 \leq i \leq l-1$, we def\/ine
\begin{gather*}
	\mathbf{m}_{Y_i^*} (n) \defeq \mathbf{m}_Y (ln +i) , \qquad n \in \Z .
\end{gather*}
Then we have an $l$-tuple
$\big( \big\{ \mathbf{m}_{Y_0^*} (n) \big\}_{n \in \Z }, \ldots ,\big\{ \mathbf{m}_{Y_{l-1}^*} (n) \big\}_{n \in \Z } \big)$ of subsequences $\big\{ \mathbf{m}_Y (\nu ) \big\}_{\nu \in \Z }$. Each $\big\{ \mathbf{m}_{Y_i^*} (n) \big\}_{n \in \Z } $ is also a Maya diagram.	
\begin{Definition}\label{def:s3.2}
The $l$-quotient for $Y$ is the $l$-tuple of Young diagrams $\underline{Y}_{(l)}^* \defeq \big( Y_0^* ,\ldots , Y_{l-1}^* \big)$ corresponding to the Maya diagrams $\big( \big\{ \mathbf{m}_{Y_0^*} (n) \big\}_{n \in \Z }, \ldots ,\big\{ \mathbf{m}_{Y_{l-1}^*} (n) \big\}_{n \in \Z } \big)$.
\end{Definition}
\begin{Remark}\label{rem:hook}
A hook of length $nl$ in $Y$ corresponds to a hook of length $n$ in $\underline{Y}_{(l)}^*$.
\end{Remark}
For the latter purpose, let us def\/ine the following
\begin{Definition}
We def\/ine $k_i (Y) \in \Z$ $(i=0, \ldots ,l-1)$ by the following condition:
\begin{gather*}%\label{eq:s10}
	\# \big\{ \mathbf{m}_{Y_i^*} (n )=1 \,|\, n <k_i (Y) \big\} =\# \big\{ \mathbf{m}_{Y_i^*} (n )=0 \,|\, n \geq k_i (Y) \big\} .
\end{gather*}
\end{Definition}
\begin{Remark}\label{rem:k_i}\quad
\begin{enumerate}\itemsep=0pt
\item[(i)] Notice that $k_i (Y) \neq 0$ in general, since we take a subsequence. However
\begin{gather*}
	\sum_{i=0}^{l-1} k_i (Y) = 0.
\end{gather*}
\item[(ii)] $ k_0 (Y) , \ldots ,k_{l-1} (Y))$ is invariant under removing a hook of length $nl$ from $Y$.
\end{enumerate}
\end{Remark}

\begin{Example}
Let us consider the $5$-quotient for the Young diagram $Y$ shown in Fig.~\ref{m}.

	(i)~In this case,
	\begin{gather*}
		\begin{array}{@{}l@{\,}llllllllllllllllllllllll}
			\{ \mathbf{m}_Y (\nu ) \} &=\{ &\underline{0} ,&0,&0,&1,&0,&1
				&| \, &0,&1,&0,&1,&1,&\underline{1} &\} ,\\
			\{ \mathbf{m}_{Y_0^*} (n) \} &=\{ &\underline{0} , &0,&&&& & \, &0,&&&& &\underline{1} &\} ,\\
			\{ \mathbf{m}_{Y_1^*} (n) \} &=\{ &\underline{0} , &&0,&&& & \, &&1,&&& &\underline{1} &\} ,\\
			\{ \mathbf{m}_{Y_2^*} (n) \} &=\{ &\underline{0} , &&&1,&& & \, &&&0,&& &\underline{1} &\} ,\\
			\{ \mathbf{m}_{Y_3^*} (n) \} &=\{ &\underline{0} , &&&&0,& & \, &&&&1,& &\underline{1} &\} ,\\
			\{ \mathbf{m}_{Y_4^*} (n) \} &=\{ &\underline{0} , &&&&&1, & \, &&&&&1, &\underline{1} &\} .
		\end{array}
	\end{gather*}
	Therefore the $5$-quotient of $Y$ is
	\begin{gather*}
		\big( Y_0^* , Y_1^* , Y_2^* , Y_3^* , Y_4^* \big)
			=\big( \varnothing \, , \, \varnothing \, , \, \Young{\Ca{}\cr} \, , \, \varnothing \, , \, \varnothing \big).
	\end{gather*}

(ii) Let us remove the hook of length $5$ from $Y$ by means of $5$-quotient $(Y_0^*,\ldots ,Y_4^*)$ of $Y$. This is equivalent to exchanging $\mathbf{m}_{Y_2^*} (-1)=1$ and $\mathbf{m}_{Y_2^*} (0)=0$, shown as follows:
\begin{gather*}
	\big\{ \mathbf{m}_{Y_2^*} (n) \big\} =\big\{ \underline{0} ,\widehat{1} \,| \, \widehat{0} ,\underline{1} \big\}
		\longmapsto \big\{ \underline{0} ,\widehat{0}\, | \, \widehat{1} ,\underline{1} \big\}.
\end{gather*}
\end{Example}

\subsubsection{Cores}
\begin{Definition}
A Young diagram is called $l$-\textit{core} if its $l$-quotient $\underline{Y}_{(l)}^*$ is empty.
We deno\-te~$\mathcal{C}^{(l)}$ by the subset of $\mathcal{P}$ which consists of $l$-cores.
\end{Definition}

\begin{Example}
The set of $2$-cores $\mathcal{C}^{(2)}$ consists of the empty partition and `stairs-like' Young diagrams which correspond to partitions $(m-1, m-2, \ldots, 1)$ for $m \geq 2$.
\end{Example}

\begin{Definition}\label{def:s3.5} Let $Y^{(l)}$ be the Young diagram obtained by removing as many hooks of length~$l$ as possible from~$Y$.\footnote{It can be shown that $Y^{(l)}$ is well-def\/ined, i.e., it does not depend on the choice of hooks to be removed.} It is called the $l$-core of~$Y$.
\end{Definition}
\begin{Proposition}\label{rem:well_core}\quad
\begin{enumerate}\itemsep=0pt
\item[$(i)$] If $k_i (Y) =0$ for all $0 \leq i \leq l-1$, then $Y^{(l)}=\varnothing$.
\item[$(ii)$] There is a bijection between $\mathcal{C}^{(l)}$ and the set
\begin{gather*}
\mathcal{K}^{(l)}\defeq \left\{ (k_{0}, \ldots, k_{l-1}) \in \Z^{l}\,\big| \, \sum_{i=0}^{l-1} k_i =0 \right\}.
\end{gather*}
\end{enumerate}
\end{Proposition}
\begin{proof}
The statement (i) follows from Remark~\ref{rem:k_i}(ii). To see~(ii), it is enough to notice that~$Y^{(l)}$ is uniquely determined by $(k_{0}(Y), \ldots, k_{l-1}(Y))$ by Remarks~\ref{rem:hook} and~\ref{rem:k_i}(ii).
\end{proof}
\begin{Proposition}[{\cite[Theorem~2.7.30]{GA}}]\label{thm:s3.6} For any $l \ge 2$, a Young diagram $Y$ is uniquely determined by its $l$-core $Y^{(l)}$ and $l$-quotient $\underline{Y}^*_{(l)}$. Thus for given $\lambda \in \mathcal{C}^{(l)}$, there is a bijection between the set
	\begin{gather*}
\mathcal{P}(n:\lambda) \defeq \bigl\{Y\in \mathcal{P}(n) \,|\, Y^{(l)} = \lambda \bigr\}
	\end{gather*}
	and $\mathcal{P}_{l}(n)$. Moreover, the weight of $Y$ is given by
	\begin{gather*}
		|Y|=\big|Y^{(l)}\big|+l\big|\underline{Y}^*_{(l)}\big| .
	\end{gather*}
\end{Proposition}

\subsubsection{Weights of cores}
Recall that we place a Young diagram $Y$ on the $(x,y)$-plane, as shown in Fig.~\ref{m}.
\begin{Definition}
For $\nu \in \Z$, we def\/ine $N_\nu (Y)$ to be the number of boxes on $y=x-\nu $.
\end{Definition}
Note that $N_{\nu}(Y)=0$ if $|\nu|$ is large enough.
\begin{Proposition}\label{prop:t1}
For $Y \in \mathcal{C}^{(l)}$ which is determined by $(k_{0}(Y), \ldots, k_{l-1}(Y))$, we have
	\begin{gather}
		|Y|= \sum_{i=0}^{l-1} \left\{ \frac{1}{2} lk_i (Y)^2 + ik_i (Y) \right\} .	\label{eq:t1}
	\end{gather}
\end{Proposition}
\begin{proof}
For $Y \in \mathcal{C}^{(l)}$, $|Y| = \sum\limits_{i=0}^{l-1} \sum\limits_{k\in \Z } N_{kl+i} (Y)$. By def\/inition, we also have
	\begin{gather*}
N_{\mu l} (Y) =
\begin{cases}
\displaystyle	\# \bigl\{ \mathbf{m}_Y (\nu ) \,|\, \mathbf{m}_Y (\nu )=1 ,\, \nu <\mu l \bigr\} \\
\displaystyle\qquad{} = \sum_i \# \bigl\{ \mathbf{m}_{Y_i^*} (n ) \,| \,\mathbf{m}_{Y_i^*} (n )=1, \, n <\mu \bigr\} , \qquad \mbox{if} \ \mu \leq 0, \\
\displaystyle \# \bigl\{ \mathbf{m}_Y (\nu ) \,|\, \mathbf{m}_Y (\nu )=0 ,\, \nu \geq \mu l \bigr\}\\
\displaystyle \qquad = \sum_i \# \bigl\{ \mathbf{m}_{Y_i^*} (n ) \,|\, \mathbf{m}_{Y_i^*} (n )=0, n \geq \mu \bigr\} , \qquad \mbox{if} \ \mu \geq 0.
		\end{cases}
	\end{gather*}
	Hence we obtain,
	\begin{gather*}
\sum_{\mu \in \Z } N_{\mu l} (Y)= N_0 (Y) + \sum_{\mu =1}^{\infty} ( N_{\mu l} (Y)+ N_{-\mu l} (Y) ) \\
\hphantom{\sum_{\mu \in \Z } N_{\mu l} (Y)}{} = \sum_{k_i >0} k_i (Y) +\sum_{n=1}^{\infty}
			\bigg\{ \sum_{k_i (Y) >n} ( k_i (Y) -n ) + \sum_{k_i (Y) <-n} ( -k_i (Y) -n ) \bigg\}.
\end{gather*}
	The right hand side of the above equation is equal to
	\begin{gather*}
\sum_{k_i (Y) > 0} \frac{1}{2} k_i (Y) ( k_i (Y) -1 ) +\sum_{k_i (Y) <0}	\frac{1}{2} ( -k_i (Y) ) ( -k_i (Y) -1 )
 +\sum_{k_i (Y) >0} k_i (Y)\\
\qquad{} =\frac{1}{2}\left( \sum_{i=0}^{l-1} k_i (Y)^2 +\sum_{i=0}^{l-1} k_i (Y) \right) =\frac{1}{2} \sum_{i=0}^{l-1} k_i (Y)^2 ,
\end{gather*}
where we have used Remark \ref{rem:k_i}(i). To compute $\sum\limits_{\mu \in \Z } N_{\mu l+n} (Y)$ for $0<n<l$, we set
\begin{gather*}
k_i^{(n)}=\begin{cases} k_i-1 & \mathrm{if} \ 0\leq i <n,\\
k_i & \mathrm{if} \ n\leq i <l.
\end{cases}
\end{gather*}
By the same argument as above, we have
\begin{gather*} \sum_{\mu \in \Z } N_{\mu l+n} (Y)= \frac{1}{2}\left( \sum_{i=0}^{l-1} k_i^{(n)} (Y)^2 +\sum_{i=0}^{l-1} k_i^{(n)} (Y) \right).
\end{gather*}
Then, using Remark \ref{rem:k_i}(i), it is easy to see that
	\begin{gather*}
		\sum_{\mu \in \Z } N_{\mu l+n} (Y)= \frac{1}{2} \sum_{i=0}^{l-1} k_i (Y)^2 + \sum_{i=n}^{l-1} k_i (Y).
	\end{gather*}
	Summing them up, we obtain the claim.
\end{proof}

\subsection[Quotients and cores for $\mathcal{P}_{\w}(\vb)$]{Quotients and cores for $\boldsymbol{\mathcal{P}_{\w}(\vb)}$}
Let $\Gamma$ be a cyclic group of order $l$. Let $\w=[W]$ be an isomorphism class of $r$-dimensional $\Gamma$-module. Let $W=\oplus _{i=1}^{r} \C \otimes \rho_{i}$ be a decomposition of $W$ into 1-dimensional $\Gamma$-modules. We regard $(\rho_1 , \ldots , \rho_r)$ as an element of $(\Z/l\Z)^r$, where we f\/ix an order of them: ${\rho}_i \leq {\rho}_{i+1}$. We give the $(\rho_1 , \ldots , \rho_r)$-coloring to $\mathcal{P}_{r}(n)$. Let us call it $\w$-\textit{coloring}. For f\/ixed $\w$, $\mathcal{P}_{r}(n)$ has the following decomposition into disjoint sets:
\begin{gather}
\mathcal{P}_{r}(n)=\bigsqcup_{\vb} \mathcal{P}_{\w}(\vb),\label{eq:disjoint}
\end{gather}
where $\mathcal{P}_{\w}(\vb)$ is the subset of $\mathcal{P}_{r}(n)$ whose number of $i$-colored boxes is equal to $v_i$ for $0\leq i \leq {l-1}$, if we give the $\w$-coloring to Young diagrams. Here $\vb = (v_0, \ldots , v_{l-1}) \in (\Z_{\geq 0})^{l}$ runs over all elements with $|\vb | \defeq \sum\limits_{i=0}^{l-1} v_i = n$.

\begin{Notation}
Hereafter in this paper, we use the following assumptions and notations on $l$-dimensional row vector.
\begin{enumerate}\itemsep=0pt
\item[(i)] We understand that indices of $l$-dimensional row vectors are in $\Z /l\Z$.
\item[(ii)] Let $\mathbf{e}_j \in \Z^l$ $(j=0,\ldots ,l-1)$ be the row vector whose $j$-th component is $1$ and the others are $0$.
\item[(iii)] Let $\boldsymbol{\delta}$ be the row vector def\/ined by $\boldsymbol{\delta}\defeq \sum\limits_{j = 0}^{l-1} \mathbf{e}_{j}$.
\end{enumerate}
\end{Notation}

Let us take $Y \in \mathcal{P}_{\mathbf{e}_j} (\vb)$. Then all boxes on the line $y=x-kl-i$ are colored by the same color $i+j$. Since $v_i$ is equal to the number of $i$-colored boxes in~$Y$, we have
\begin{gather}
	v_i =\sum_{\mu \in \Z } N_{\mu l+i-j} (Y).\label{eq:v}
\end{gather}

\begin{Proposition}
	All Young diagrams in $\mathcal{P}_{\mathbf{e}_j} (\vb )$
	have the same $l$-core.
\label{prop:s3.8}
\end{Proposition}	
\begin{proof}
For $Y \in \mathcal{P}_{\mathbf{e}_j} (\vb )$,	let $\left\{ \mathbf{m}_Y (\nu ) \right\}_{\nu =-\infty}^{\infty} $ be the Maya diagram corresponding to~$Y$. For all $m \leq 0$, we have
	\begin{gather*}
N_{m} (Y)= \# \bigl\{ \mathbf{m}_Y (\nu ) \,|\, \mathbf{m}_Y (\nu )=1 ,\, \nu <m \bigr\} ,
	\end{gather*}
	and for all $m \geq 0$, we have
	\begin{gather*}
		N_{m} (Y)= \# \bigl\{ \mathbf{m}_Y (\nu ) \,|\, \mathbf{m}_Y (\nu )=0 ,\,\nu \geq m \bigr\} .
	\end{gather*}
	Then we have
	\begin{gather*}
v_{i+j}-v_{i+j+1} = \sum_{\mu \in \Z} N_{\mu l+i} (Y)- \sum_{k \in \Z} N_{\mu l+i+1} (Y) = \sum_{\mu \in \Z} ( N_{\mu l+i} (Y) - N_{\mu l+i+1} (Y)) \\
\hphantom{v_{i+j}-v_{i+j+1}}{} = \sum_{\mu \geq 0} ( N_{\mu l+i} (Y) - N_{\mu l+i+1} (Y) ) -\sum_{\mu < 0} (N_{\mu l+i+1} (Y) - N_{\mu l+i} (Y) ) \\
\hphantom{v_{i+j}-v_{i+j+1}}{} = \sum_{\mu \leq 0} (1-\mathbf{m}_Y (\mu l +i) ) - \sum_{\mu <0} \mathbf{m}_Y (\mu l+i ) \\
\hphantom{v_{i+j}-v_{i+j+1}}{} = \# \bigl\{ \mathbf{m}_Y ( \mu l+i ) =0 \,|\, \mu \in \Z_{\geq 0}\bigr\}
						- \# \bigl\{ \mathbf{m}_Y (\mu l+i ) =1 \,|\, \mu \in \Z_{< 0} \bigr\} \\	
\hphantom{v_{i+j}-v_{i+j+1}}{} = \# \bigl\{ \mathbf{m}_{Y_i^*} (\mu ) =0 \,|\, \mu \in \Z_{\geq 0} \bigr\}
						- \# \bigl\{ \mathbf{m}_{Y_i^*} (\mu ) =1\,|\, \mu \in \Z_{< 0} \bigr\} =k_{i} (Y).
	\end{gather*}
Now our claim follows from Proposition~\ref{rem:well_core}(ii).
\end{proof}

By def\/inition, a hook of length $nl$ has $n$ $i$-colored boxes for all $0\leq i\leq l-1$. We have the following
\begin{Corollary}\label{cor:QforC}
For each $0\leq j\leq l-1$,
\begin{gather*}
\mathcal{P}_{\mathbf{e}_j} (n\boldsymbol{\delta}) =\mathcal{P} (nl:\varnothing ).
\end{gather*}
In particular $\mathcal{P}_{\mathbf{e}_j} (n\boldsymbol{\delta})$ is independent of~$j$.
\end{Corollary}
\begin{proof}\looseness=1 In the proof of Proposition \ref{prop:s3.8}, we obtained the formula $k_i (Y) =v_{i+j} -v_{i+j+1}$. Thus, for $Y \in \mathcal{P}_{\mathbf{e}_j} (n\boldsymbol{\delta})$, $k_i (Y)=0$ for all $0\leq i \leq l-1$. This means that $Y^{(l)} =\varnothing$, by Proposi\-tion~\ref{rem:well_core}(i).
\end{proof}

\section{Quiver varieties and torus actions}\label{sec:s2}

In this section, we introduce quiver varieties and torus actions on them. We follow~\cite{Nakajima3} and~\cite{VV} closely.

\subsection{Quiver varieties}

\subsubsection{Framed moduli space of torsion free sheaves}
 Let $V$, $W$ be complex vector spaces with $\dim V=n$, $\dim W=r$. We consider
 \begin{gather*}
 \mathbf{M}_r (n) \defeq (Q \otimes \Hom (V,V))\oplus \Hom (W, V) \oplus \Hom (V,W),
 \end{gather*}
 where $Q$ is a 2-dimensional complex vector space. We def\/ine an action of $\GL (V)$ on $\mathbf{M}_{r} (n)$ by
 \begin{gather*}
 (B, i, j) \mapsto \big((\mathrm{id}_Q \otimes g)B(\mathrm{id}_Q \otimes g)^{-1}, ig^{-1}, gj\big), \qquad g \in \GL (V).
 \end{gather*}
 We def\/ine a map $\mu\colon \mathbf{M}_{r} (n) \to \End (V)$ by
 \begin{gather*}
 \mu(B, i, j) \defeq [B\wedge B] + ij .
 \end{gather*}
 Note that $\mu^{-1}(0)$ is invariant under the action of $\GL (V)$.
 \begin{Definition} \label{def:stability}
 We say $(B, i, j) \in \mathbf{M}_{r} (n)$ is \textit{stable} if the following condition is satisf\/ied:
 If a subspace $S$ of $V$ contains $\mathrm{Im}(i) $ and $B(S) \subset Q \otimes V$, then $S=V$.
 \end{Definition}

We def\/ine
 \begin{gather*}
 \M_{r} (n) \defeq \big\{(B,i,j)\in \mu^{-1}(0) \,|\, (B, i, j) \ \mathrm{is} \ \mathrm{stable} \big\} / \GL (V) ,
 \end{gather*}
 and
 \begin{gather*}
 \M_{r}(n)_0 \defeq \mu^{-1}(0) /\negthickspace/ \GL(V)= \mbox{the set of \it{closed}} \ \GL(V)\mbox{-orbits in} \ \mu^{-1}(0).
\end{gather*}
 These are \textit{projective quotient} and \textit{affine quotient} respectively, in the sense of geometric invariant theory (GIT). By the general theory of GIT, There is a proper morphism $\pi \colon \M_{r} (n) \to \M_{r}(n)_0$. The following theorem is proved in \cite[Section~2]{NY}.

\begin{Theorem} \label{thm:NY}\quad
\begin{enumerate}\itemsep=0pt
\item[$(i)$] $\M_{r} (n)$ is a nonsingular complex algebraic variety of dimension $2nr$ and $\pi \colon \M_{r} (n) \to \M_{r}(n)_0$ is a resolution of singularities.
\item[$(ii)$] $\M_{r} (n) $ is isomorphic to the moduli space of the pair $(E, \phi)$ where $E$ is a torsion free sheaf over $\CP^2$ of $\rank E = r$,
 $c_2 (E) = n$ which is locally free in a neighborhood of $l_{\infty} = \{[0:z_{1}: z_{2}]\}$ and $\phi$ is an isomorphism $E|_{l_{\infty}} \to
\shfO_{\CP^2}^{\oplus r}$ $($framing at infinity$)$.
\end{enumerate}
 \end{Theorem}
 \begin{Remark}
When $\dim (W) =1$, any torsion free sheaf $E$ of rank 1 with the above condition is a subsheaf of $\shfO_{\CP^2}$ such that $\Supp(\shfO_{\CP^2}/E)$ is a $0$-dimensional subscheme of $\C^2$, so we recover the Hilbert schemes of points on $\C^2\colon \Hilbn{\big(\C^2\big)} \cong \M_{1} (n)$. Therefore, $\M_{r} (n)$ is a higher rank generalization of the Hilbert scheme of points. Note that $\M_{r}(n)_0$ is isomorphic to the $n$-th symmetric product $S^n\big(\C^2\big)$ of $\C^2$ and the above morphism $\pi$ coincides with the Hilbert--Chow morphism $\pi \colon \Hilbn{\big(\C^2\big)} \to S^n\big(\C^2\big)$ (see Appendix \ref{sub:app}).
\end{Remark}

\subsubsection{Quiver varieties}\label{sub:quiver}
Let $\Gamma$ be a f\/inite cyclic subgroup of $\SL (2,\C )$ of order $l$. Let $Q$ be the 2-dimensional $\Gamma$-module def\/ined by the inclusion $\Gamma \subset \SL (2, \C)$. For $\Gamma$-modules $V$ and $W$, let
 \begin{gather*}
 \mathbf{M}_{W} (V) \defeq (Q \otimes \Hom (V,V))_{\Gamma} \oplus \Hom_{\Gamma}(W,V) \oplus \Hom_{\Gamma}(V,W),
 \end{gather*}
where $(\; \;)_{\Gamma}$ means $\Gamma$-invariant part. We def\/ine $\GL _{\Gamma}(V)$-action on $\mathbf{M}_{W} (V)$, and the map $\mu \colon \mathbf{M}_{W} (V) \to \End _{\Gamma}(V)$ as above. We also def\/ine the stability condition by the same condition as Def\/inition \ref{def:stability}, where subspace $S \subset V$ is replaced by $\Gamma$-submodule. Then we def\/ine $\M_{\w} (\vb)$ and $\M_{\w} (\vb)_\mathbf{0}$ exactly the same as above, where $\vb$ and $\w$ are isomorphism classes of $V$ and $W$ as $\Gamma$-module respectively. Let $\{R_i\}_{i=0}^{l-1}$ be the set of all irreducible representations of $\Gamma$,
where $R_0$ is the trivial representation. We can denote $V=\bigoplus_{i=0}^{l-1} V_i \otimes R_i$, $W=\bigoplus_{i=0}^{l-1} W_i \otimes R_i$. Then $\vb$ and $\w$ can be regarded as row vectors $( \dim V_0 ,\ldots ,\dim V_{l-1} )$ and $( \dim W_0 ,\ldots ,\dim W_{l-1} )$ respectively. We call $\M_{\w} (\vb)$ and $\M_{\w} (\vb)_\mathbf{0}$ \textit{quiver varieties}. The restriction of $\pi \colon \M_{r}(n) \to \M_{r}(n)_{0}$ gives a~proper morphism $\pi \colon \M_{\w} (\vb) \to \M_{\w} (\vb)_\mathbf{0}$.
\begin{Theorem}\label{thm:M_w(vb)}\quad
\begin{enumerate}\itemsep=0pt
\item[$(i)$] $\M_{\w} (\vb)$ is a nonsingular complex algebraic variety and $\pi \colon \M_{\w} (\vb) \to \M_{\w} (\vb)_\mathbf{0}$ is a resolution of singularities. If $\M_{\w} (\vb) \neq \varnothing$, its dimension is given by
 \begin{gather*}
	\dim_{\C} \M_{\w} (\vb )= \vb \widetilde{\mathbf{C}}_{\Gamma} {}^t \vb +2\vb^t \w ,	%\label{eq:s4}
 \end{gather*}
 where $\widetilde{\mathbf{C}}_{\Gamma}$ is the affine Cartan matrix of type $A_{l-1}^{(1)}$
 and ${}^t (\, )$ means the transposition.
\item[$(ii)$] $\M_{\w} (\vb)$ is isomorphic to the moduli space parameterizing the pair~$(E, \phi)$, where $E$ is a~$\Gamma$-equivariant torsion free sheaf with $H^1(\CP^2, E\otimes \shfO_{\CP^2}(-l_\infty)) \cong V$ $($as a~$\Gamma$-module$)$ and~$\phi$ is a~$\Gamma$-equivariant trivialization $E|_{l_\infty} \to \shfO_{l_{\infty}} \otimes W$.
\end{enumerate}
\end{Theorem}

The statements (i) is due to Nakajima \cite{Nakajima2}. The statement (ii) is stated in \cite[Section~2.3]{VV}. Note that there is a natural $\Gamma$-action on $\M_{r} (n)$ induced by that on~$\CP^2$, if we f\/ix a lift of the $\Gamma$-action to $\shfO_{l_{\infty}}^{\oplus r}$ by~$\w$. Then the above theorem tells us that $\Gamma$-f\/ixed point set $\M_{r} (n)^{\Gamma}$ of $\M_{r} (n)$ have the following decomposition:
\begin{gather*}
	\M_{r} (n)^{\Gamma} = \bigsqcup_{\vb \atop |\vb|=n} \M_{\w} (\vb).%\label{eq:decomp}
\end{gather*}

\subsubsection{}\label{sub:ale}
For the latter purpose, let us introduce the so-called \textit{ALE spaces}.
\begin{Theorem}[{\cite{IN1, Kronheimer}}]\label{rem:ale}\quad
\begin{enumerate}\itemsep=0pt
\item[$(i)$] $\M_{\mathbf{e}_0} (\boldsymbol\delta)_{\mathbf{0}}$ is isomorphic to the simple singularity ${\C}^2 / \Gamma$.
\item[$(ii)$] $\pi\colon \M_{\mathbf{e}_0}(\boldsymbol\delta) \to \M_{\mathbf{e}_0} (\boldsymbol\delta)_{\mathbf{0}}$ is the minimal resolution. Here we regard $\boldsymbol\delta$ as the isomorphism class of the regular representation of $\Gamma$. $\M_{\mathbf{e}_0}(\boldsymbol\delta)$ is called \textit{ALE space}.\footnote{ALE represents \textit{asymptotically locally Euclidean}, since $\M_{\mathbf{e}_0}(\boldsymbol\delta)$ admits a Riemannian metric which approximates the f\/lat metric on $\C^2/\Gamma$ at inf\/inity.}
\end{enumerate}
\end{Theorem}
This result was f\/irst obtained by Kronheimer \cite{Kronheimer}, and rediscovered by Ito--Nakamura \cite{IN1} and Ginzburg--Kapranov (unpublished) independently.

\begin{Remark}
Results in Sections~\ref{sub:quiver} and~\ref{sub:ale} holds for any f\/inite subgroup of $\SL (2,\C )$, with a minor modif\/ication.
\end{Remark}

\subsection{Torus actions and their f\/ixed point set}

\subsubsection{Torus actions on framed moduli spaces}
Let $T_Q$ (resp.~$T_W$) be a maximal torus in $\GL(Q)$ (resp.~$\GL(W)$), and let $T = T_Q \times T_W$. Note that $T_Q \cong (\C^{*})^2$ and $T_W \cong (\C^{*})^r$. Let us consider the following $T$-action on $\mathbf{M}_{r} (n)$:
\begin{gather*}
(B_1, B_2, i, j) \mapsto \big(t_{1} B_{1}, t_{2} B_{2}, ie^{-1}, t_{1}t_{2}ej\big) ,
\end{gather*}
where $t=(t_{1}, t_{2})\in T_Q$, $e=\operatorname{diag}(e_1, \ldots, e_{r})\in T_W$, and we denote $B=(B_1, B_2)\in Q \otimes \Hom (V, V)$. This induces a $T$-action on $\M_{r} (n)$ and $\M_{r}(n)_{0}$. Note that the morphism $\pi \colon \M_{r} (n)$ $\to \M_{r}(n)_{0}$ is $T$-equivariant with respect to this action. This torus action is studied in~\cite{NY}. It turns out that $T$-f\/ixed points in~$\M_{r}(n)$ are parameterized by $\mathcal{P}_{r}(n)$ and $T$-f\/ixed points in $\M_{r}(n)_0$ consists of single point.

\subsubsection{Torus actions on quiver varieties}
We can regard $T_Q$ (resp.~$T_W$) as a maximal torus in $\GL_{\Gamma}(Q)$ (resp.~$\GL_{\Gamma}(W)$). Thus we see that $T = T_Q \times T_W$ also acts on $\mathbf{M}_{W}(V)$ and this induces an action of $T$ on quiver varieties $\M_{\w} (\vb)$ and $\M_{\w}(\vb)_{\mathbf{0}}$. We can study $T$-f\/ixed points on quiver varieties by using results in~\cite{NY}.
\begin{Proposition}\label{prop:t6}\quad
\begin{enumerate}\itemsep=0pt
\item[$(i)$] $(E, \phi) \in \M_{r} (n)^{\Gamma}$ is fixed by the $T$-action if and only if $E$ has a decomposition $E = I_1 \oplus \cdots \oplus I_r$ satisfying the following conditions for $\alpha = 1, \ldots , r$.
\begin{enumerate}\itemsep=0pt
\item[$(a)$] $I_\alpha$ is an ideal sheaf of $0$-dimensional subscheme $Z_\alpha$ whose support is a $\Gamma$-orbits in $\C^2$.
\item[$(b)$] Under $\phi$, $I_\alpha |_{l_\infty}$ is mapped to the $\alpha$-th factor of $\shfO_{l_\infty} \otimes W$.
\item[$(c)$] $I_\alpha$ is fixed by the $T_Q$-action.
\end{enumerate}
\item[$(ii)$] The $T$-fixed point set in $\M_{\w} (\vb)_{\mathbf{0}}$ consists of a single point.
\end{enumerate}
\end{Proposition}
\begin{proof}
This follows from the same argument as \cite[Proposition 2.9]{NY}.
\end{proof}

\subsubsection{Parameterization of f\/ixed points}
It is instructive to compare equation~(\ref{eq:disjoint}) with equation~(\ref{eq:decomp}).
\begin{Proposition}\label{prop:fixed}
There is a one-to-one correspondence between the $T$-fixed point set of $\M_{\w} (\vb)$ and $\mathcal{P}_{\w}(\vb)$.
\end{Proposition}

\begin{proof} This follows from Proposition \ref{prop:t6} and the same argument as in \cite[Proposition~2.9(2)]{NY}. The $T$-f\/ixed points $(E=\oplus_{\alpha} I_{\alpha}, \phi) \in \M_{\w} (\vb)$ are parametrized by $\mathcal{P}_{\w}(\vb)$ as follows. By the conditions~(a) and~(c) in Proposition~\ref{prop:t6}(i), $I_{\alpha}$ is identif\/ied with an ideal in $\mathbb{C}[x,y]$ generated by monomials. Hence it corresponds to a~Young diagram $Y_{\alpha}$ as in \cite[Proposition~2.9(2)]{NY}. Note that $\mathbb{C}[x,y]/I_{\alpha}$ is a $\Gamma$-module via the $2$-dimensional representation~$Q$ of~$\Gamma$. It is standard that $V=H^1(\mathbb{CP}^2, E(-\ell_{\infty}))$ is isomorphic to $\oplus_{\alpha} \mathbb{C}[x,y]/I_{\alpha}$ as a~$\mathbb{C}$-vector space. (This follows from Lemma~2.2 in \cite[Chapter~2]{Lecture}.) Via the ADHM description, $\Gamma$-equivariant framing $\phi$ in
Proposition \ref{prop:t6}(i)(b) corresponds to a $\Gamma$-equivariant linear map $i\colon W\to V$. (See the proof of Theorem~4.4 in \cite[Chapter~4]{Lecture}.) It then follows that $V \cong \oplus_{\alpha} ( \mathbb{C}[x,y]/I_{\alpha} ) \otimes \rho_{\alpha}$ as a~$\Gamma$-module, where $W=\oplus_{\alpha} \mathbb{C}\otimes \rho_{\alpha}$. Hence if we color the $r$-tuple of Young diagrams $\underline{Y}=(Y_{\alpha})$ by $\w=(\rho_{\alpha})$, it follows that $\underline{Y}\in \mathcal{P}_{\w}(\vb)$ since the isomorphism class of the $\Gamma$-module $\oplus_{\alpha} ( \mathbb{C}[x,y]/I_{\alpha} ) \otimes \rho_{\alpha}$ is $\vb=[V]$.
\end{proof}

\subsubsection{}
We compare Proposition \ref{prop:fixed} with the work of Nakajima \cite{Nakajima1}. The results are not used later. In~\cite{Nakajima1}, an action of a certain 1-parameter subgroup of $T$ on $\M_{\w} (\vb )$ is considered. And in the case of $v_0 = w_0 =0$, f\/ixed points are parametrized by Young tableaux.

\begin{Definition} Let $\mu = (\mu_1, \ldots, \mu_{l})$ be an $l$-tuple of nonnegative integers.
\begin{enumerate}\itemsep=0pt
\item[(i)] A $\mu$\textit{-tableau of shape} $Y$ is a Young diagram $Y$ whose boxes are numbered with the f\/igures from $1$ to $l$ such that the cardinality of the boxes with f\/igure $k$ is $\mu_k$.
\item[(ii)] A $\mu$-tableau of shape $Y$ is said to be {\it row-increasing} if the entries in each row increase strictly from the left to the right.
\end{enumerate}
\end{Definition}
\begin{Proposition}[{\cite[Lemma 5.8]{Nakajima1}}]
When $v_0 = w_0 =0$, $T$-fixed points on $\M_{\w} (\vb)$ are parametrized by row increasing $\mu$-tableaux of shape $\lambda =(1^{w_1} ,2^{w_2} ,\ldots ,(l-1)^{w_{l-1}} )$,\footnote{This means that $\lambda$ is a partition with $w_i = \# \{ j \,|\, \lambda_j = i \}$.} where ${\mu}_k \defeq v_n +\sum\limits_{i\geq k} u_i \, (k=1,\ldots ,l-1)$ and ${\mu}_{l} \defeq v_{l-1}$.
\end{Proposition}
\begin{Proposition}
When $v_0 = w_0 =0$, there is a one-to-one correspondence between $\mathcal{P}_{\w} (\vb)$ and the set of row increasing $\mu$-tableau of shape $\lambda= (\lambda_1,\ldots , \lambda_r)\defeq(1^{w_1} ,2^{w_2} ,\ldots ,(l-1)^{w_{l-1}} )$.
\end{Proposition}
\begin{proof}
Note that, in this case, $\mathcal{P}_{\w} (\vb)$ consists of tuples of $l$-cores. Each row in a $\mu$-tableau of shape $\lambda$ corresponds to a $l$-core.
A correspondence is given as follows. For the $i$-th row of a~$\mu$-tableau of shape $\lambda$, we denote $\mathfrak{C}_i$ by the set of contents on it. The data $\left(k_{0}, \ldots, k_{l-1}\right)$ of a~$l$-core is given as follows. For $1 \leq n \leq \lambda_i$,
\begin{gather*}
	k_{n-\lambda_i-1} =
	\begin{cases}
		\hphantom{-}0, & n \in \mathfrak{C}_i , \\
		-1 & n \not\in \mathfrak{C}_i,
	\end{cases}
\end{gather*}
and for $\lambda_i < n \leq l$,
\begin{gather*}
	k_{n-\lambda_i-1} =
	\begin{cases}
		1 & n \in \mathfrak{C}_i, \\
		0 & n \not\in \mathfrak{C}_i.
	\end{cases}
\end{gather*}
Then it is easy to see that this is a bijection.
\end{proof}

\subsubsection{Characters at f\/ixed points}\label{sec:s4}
\begin{Notation}\quad
\begin{enumerate}\itemsep=0pt
\item[(i)] We denote by $e_1 ,\ldots ,e_r $ $1$-dimensional $T$-modules given by $T \to \C^{*}$, $t =(t_1, t_2, e_1, \ldots, e_r )$ $\mapsto e_{i}$. Similar notations are used for $t_1$ and $t_2$.
\item[(ii)] For $a,b \in \Z$, we set
\begin{gather*}
	\delta^{(l)}_{a,b} \defeq
	\begin{cases}
		1 & \mbox{if} \ a \equiv b \ (\mbox{mod} \ l), \\
		0 & \mbox{if} \ a \not\equiv b \ (\mbox{mod} \ l).
	\end{cases}
\end{gather*}
\end{enumerate}
\end{Notation}
Let $W=\oplus _{i=1}^{r} \C \otimes \rho_{i}$ be a decomposition of a representative $W$ of $\w$ into 1-dimensional $\Gamma$-modules as in Section~\ref{subsec:2.4}. Then we have the following
\begin{Theorem}\label{thm:s4.3}
Let $x$ be a $T$-fixed point of $\M_{\w} (\vb)$ corresponding to $( Y_1, \ldots, Y_r) \in \mathcal{P}_{\w} (\vb )$. Then the $T$-module structure of $T_x \M_{\w} (\vb)$	is given by
	\begin{gather*}
		T_x \M_{\w} (\vb)= \sum_{\alpha,\beta =1}^r N_{e_\alpha ,e_\beta }^{\Gamma } (t_1 ,t_2 ),	%\label{eq:s18}
	\end{gather*}
where
	\begin{gather*}
N_{e_\alpha ,e_\beta }^{\Gamma } (t_1 ,t_2 )= e_\beta e_\alpha^{-1}\Bigg\{ \sum_{s \in Y_\alpha }\Big( t_1^{-l_{Y_{\beta}} (s)} t_2^{a_{Y_\alpha} (s) +1} \Big)\delta^{(l)}_{h_{Y_{\beta},Y_\alpha} (s) , \rho_\alpha -\rho_\beta } \\
\hphantom{N_{e_\alpha ,e_\beta }^{\Gamma } (t_1 ,t_2 )	= e_\beta e_\alpha^{-1} \Bigg\{ }{} + \sum_{t \in Y_{\beta}}
				\Big( t_1^{l_{Y_\alpha} (t)+1} t_2^{-a_{Y_{\beta}} (t) } \Big)
					\delta^{(l)}_{h_{Y_\alpha,Y_{\beta}} (t) , \rho_\beta -\rho_\alpha }\Bigg\} .%	\label{eq:s19}
	\end{gather*}
\end{Theorem}

\begin{proof} This follows by taking $\Gamma$-invariant part of the $T$-module given in \cite[Theorem~2.11]{NY}.
\end{proof}

\section{Enumerative geometry on quiver varieties} \label{sec:enu}
This section is devoted to computations of some global invariants of quiver varieties. In Section~\ref{sec:integral}, we compute the `equivariant volume' of $\M_{\w} (\vb)$. This is a simple application of Theo\-rem~\ref{thm:s4.3} and localization theorem. The main part in this section is Section~\ref{sub:euler}, where we consider generating functions of Poincar\'e polynomials and Euler characteristics of quiver va\-rie\-ties. Our computation is based on combinatorial arguments over \textit{cores} and \textit{quotients} of Young diagrams in Section~\ref{sec:s3}.

\subsection[Equivariant integrals on $\M_{\w} (\vb)$]{Equivariant integrals on $\boldsymbol{\M_{\w} (\vb)}$}\label{sec:integral}

\subsubsection{Localization}\label{sub:vol}
First, we review some basic facts about equivariant integrals and the localization formula. Let~$M$ be an algebraic variety with an action of a torus $T$ and $H^{T}_{*}(M)$ be the equivariant homology\footnote{See \cite[Appendix~C]{NY2} for a precise def\/inition.}(with $\Q$-coef\/f\/icient). The ring $H^{T}_{*}(M)$ is a module over $H^{T}_{*}(pt)$. We consider the $T$-equivariant integral $\int_{M} 1$. This integral takes values in the quotient f\/ield of~$H^{*}_{T}(pt)$. We refer the reader \cite[Section~4.1]{NY2} for a precise def\/inition of the integral. Assume that $M$ is smooth and $T$-f\/ixed point set $M^{T}$ is f\/inite. Then the classical Atiyah--Bott localization theorem \cite{AB} says that
\begin{gather*}
\int_{M} 1 = \sum_{x \in M^T} \frac{1}{e_x},
\end{gather*}
where $e_x$ is the $T$-equivariant Euler class of $T_{x}M$. ($e_x \neq 0$ since $x$ is isolated.)

\subsubsection{}
We apply the localization formula to $\M_{\w}(\vb)$. In this case, the integral is def\/ined by equivariant pushforward of the fundamental class $[\M_{\w}(\vb)]$ to the unique $T$-f\/ixed point $o$ in $\M_{\w}(\vb)_\mathbf{0}$ by $\pi \colon \M_{\w}(\vb) \to \M_{\w}(\vb)_\mathbf{0}$ (see Proposition~\ref{prop:t6}(ii)). Furthermore, equivariant Euler classes at f\/ixed points are given by Theorem~\ref{thm:s4.3}. Let $({\epsilon}_1 ,{\epsilon}_2 , \vec{a} )$ be a coordinate on the Lie algebra of $T$ given by $t_1 =e^{{\epsilon}_1}$, $t_2 =e^{{\epsilon}_2}$, $e_\alpha =e^{a_\alpha}$, and $\vec{a} =(a_1 ,\ldots ,a_r )$.
\begin{Proposition}\label{cor:s4.4}	The $T$-equivariant integral on $\M_{\w} (\vb)$ is given by
\begin{gather*}
\int_{\M_{\w} (\vb)} 1 =\sum_{\{Y_\alpha \}_{\alpha=1}^r \in \mathcal{P}_{\w} (\vb )} \frac{1}{\prod\limits_{\alpha,\beta=1}^r n^{\Gamma}_{e_\alpha e_\beta} ({\epsilon}_1 ,{\epsilon}_2 ,\vec{a} )} ,
	\end{gather*}
where
\begin{gather*}
		n^{\Gamma}_{e_\alpha e_\beta} ({\epsilon}_1 ,{\epsilon}_2 ,\vec{a} ) =
			\prod_{s \in Y_\alpha \atop h_{Y_\beta Y_\alpha} (s) \equiv \rho_\alpha -\rho_\beta \, {\rm{mod}} \, l}
			 ( -l_{Y_\beta} (s) {\epsilon}_1 + ( a_{Y_\alpha} (s)+1 ) {\epsilon}_2 +a_\beta -a_\alpha ) \\
\hphantom{n^{\Gamma}_{e_\alpha e_\beta} ({\epsilon}_1 ,{\epsilon}_2 ,\vec{a} ) =}{} \times
			\prod_{t \in Y_\beta \atop h_{Y_\alpha Y_\beta} (t) \equiv \rho_\beta -\rho_\alpha \, {\rm{mod}} \, l}
			 ( ( l_{Y_\alpha} (t)+1 ) {\epsilon}_1 -a_{Y_\beta} (t){\epsilon}_2 + a_\beta - a_\alpha ).
	\end{gather*}
\end{Proposition}
\begin{Example}\label{ex:s4.8}
In rank $1$ cases, we have
\begin{gather*}
\int_{\mathcal{M}_{\mathbf{e}_{j}}(\vb)} 1 = \frac{1}{n! l^n ({\epsilon}_1 {\epsilon}_2 )^n} ,
\end{gather*}
where $n = |\underline{Y}_{(l)}^*|$ for $Y \in \mathcal{P}_{\mathbf{e}_j}(\vb)$. A proof of this identity is given in Proposition~\ref{prop:s4.6}.
\end{Example}
\begin{Example} Let us take $r=2$, $l=2$, $\vb =(1,1)$, and $\w =(1,1)$. There are f\/ive elements in~$\mathcal{P}_{\w}(\vb)$, as shown in Table~\ref{eq:s6}.
\begin{table}[t]\centering
	\begin{tabular}{|c||c|c|c|c|c|} \hline
	& I & II & III & IV & V \\ \hline
	$Y_1$ & $\Young{\Box{}&\Box{}\cr}$ & $\Young{\Box{}\cr \Box{}\cr}$ & ${\varnothing}$ & ${\varnothing}$ & 	$\Young{\Box{}\cr}$ \\ \hline
	$Y_2$ & ${\varnothing}$ & ${\varnothing}$ & $\Young{\Box{}&\Box{}\cr}$ & $\Young{\Box{}\cr \Box{}\cr}$ & $\Young{\Box{}\cr}$ \\ \hline
	\end{tabular}
	\caption{The elements of $\mathcal{P}_{\w}(\vb)$ for $\w=(1,1)$ and $\vb =(1,1)$.}	\label{eq:s6}
\end{table}
Characters at these f\/ixed points are given as follows:
	\begin{alignat*}{3}
			&\mbox{I}) \quad & & t_1 t_2^{-1} +t_2^2 + e_1 e_2^{-1} t_2^{-1} +e_2 e_1^{-1} t_1 t_2^2, & \\
			&\mbox{II}) \quad & & t_1^{-1} t_2 + t_1^2 + e_2 e_1^{-1} t_1^2 t_2 + e_1 e_2^{-1} t_1^{-1}, & \\
			&\mbox{III}) \quad & & t_2^2 + t_1 t_2^{-1} +e_2 e_1^{-1} t_2^{-1} +e_1 e_2^{-1} t_1 t_2^2, & \\
			&\mbox{IV}) \quad & & t_1^{-1} t_2 +t_1^2 + e_2 e_1^{-1} t_1^{-1} +e_1 e_2^{-1} t_1^2 t_2, & \\
			&\mbox{V}) \quad & & e_1 e_2^{-1} t_2 + e_2 e_1^{-1} t_1 +e_2 e_1^{-1} t_2 +e_1 e_2^{-1} t_{1}, &
		\end{alignat*}
and we have
\begin{gather*}
	\int_{\M_{\w}(\vb)} 1 = \frac{4{\epsilon}_1^2 +10{\epsilon}_1 {\epsilon}_2 +4{\epsilon}_2^2 -(a_1 -a_2 )^2}
		{{\epsilon}_1 {\epsilon}_2 E(1,-1,2,1) E(1,-1,1,2) E(-1,1,2,1) E(-1,1,1,2)} ,
\end{gather*}
where $E(x,y,z,w)\defeq xa_1 +ya_2 +z{\epsilon}_1 +w{\epsilon}_2$.
\end{Example}
\begin{Remark}
It seems interesting to study the generating function of equivariant integrals (see Sectiosn~\ref{sub:n=2} and~~\ref{sec:example}).
\end{Remark}

\subsection[Poincar\'e polynomials and Euler characteristics of $\M_{\w} (\vb)$]{Poincar\'e polynomials and Euler characteristics of $\boldsymbol{\M_{\w} (\vb)}$}\label{sub:euler}

\subsubsection{Poincar\'e polynomials for rank 1 case}
First, we consider Poincar\'e polynomials for rank $1$ case. Let us consider the following generating function of Poincar\'e polynomials:
\begin{gather*}
	\mathcal{Z}_{\mathbf{e}_{j}} (\mathfrak{t},\mathfrak{q},\vec{\mathfrak{r}} ) \defeq \sum_{m\geq 0 \atop \vb \in (\Z_{\geq 0})^{l}}
	\mathrm{b}_m ( \M_{\mathbf{e}_j}(\vb))\mathfrak{t}^m \mathfrak{q}^{|\vb|} \prod_{i=0}^{l-1} \mathfrak{r}_{i}^{v_i} ,
\end{gather*}
where $\vec{\mathfrak{r}}=(\mathfrak{r}_{0}, \ldots, \mathfrak{r}_{l-1})$ and $\mathrm{b}_m$ is the $m$-th Betti number. By the same argument as \cite[Section~3.3]{NY2}, we have
\begin{gather}
	\mathcal{Z}_{\mathbf{e}_{j}} (\mathfrak{t},\mathfrak{q},\vec{\mathfrak{r}}) = \sum_{\vb \in (\Z_{\geq 0})^{l}} \mathfrak{q}^{|\vb|}
	\prod_{i=0}^{l-1} \mathfrak{r}_{i}^{v_i} \sum_{ Y \in \mathcal{P}_{\mathbf{e}_j} (\vb) } \mathfrak{t}^{ \# \{ s \in Y \,|\, h(s) \equiv 0\, (\mathrm{mod} \, l ), \, l(s)>0 \} }.
\label{eq:t4}
\end{gather}
\begin{Example}
When $l=3$, $\vb =(1,2,1)$, $\w =(0,1,0)$,
\begin{gather*}
	\mathcal{P}_{\w} (\vb )= \left\{
	\begin{matrix}
	\Young{
	 \Box{}&\Box{\star }&\Box{}&\Box{}&\cr}, \,
	\Young{
 \Box{}&\Box{}&\cr
 \Box{\star }&\Box{}&\cr}, \,
	\Young{
	 \Box{}&\cr
 \Box{}&\cr
 \Box{\diamond }&\cr
	 \Box{}&\cr}
	\end{matrix}
	\right\}.
\end{gather*}
In the right hand of (\ref{eq:t4}), $\star$'s are counted but $\diamond$ is not.
\end{Example}
From the combinatorial discussions in Section~\ref{sec:s3}, we have the following factorization formula for $\mathcal{Z}_{\mathbf{e}_{j}} (\mathfrak{t},\mathfrak{q},\vec{\mathfrak{r}})$.
\begin{Theorem}\label{thm:t2}
	\begin{gather*}
	\mathcal{Z}_{\mathbf{e}_{j}} (\mathfrak{t},\mathfrak{q},\vec{\mathfrak{r}}) =\mathcal{Z}^{\rm{quot}} (\mathfrak{t},\mathfrak{q},\vec{\mathfrak{r}}) \, \mathcal{Z}_{\mathbf{e}_{j}}^{\rm{core}} (\mathfrak{q},\vec{\mathfrak{r}}),
	\end{gather*}
	where
\begin{gather*}
	 \mathcal{Z}^{\rm{quot}} (\mathfrak{t},\mathfrak{q},\vec{\mathfrak{r}})
		\defeq \sum_{\vb \in {\Z}^{l}_{\geq 0} \atop v_0= \cdots = v_{l-1} } \mathfrak{q}^{|\vb|}
		\prod_{i=0}^{l-1} \mathfrak{r}_{i}^{v_i}
		\sum_{ Y \in \mathcal{P}_{\mathbf{e}_j }(\vb) }
		\mathfrak{t}^{ \# \{ s \in Y \,|\, h(s) \equiv 0\, ({\rm{mod}} \, l ),\, l(s)>0 \} },
	\end{gather*}
and
\begin{gather*}
	\mathcal{Z}^{\rm{core}}_{\mathbf{e}_j} (\mathfrak{q},\vec{\mathfrak{r}}) \defeq \sum_{Y \in \mathcal{C}^{(l)}} \mathfrak{q}^{|Y|} \prod_{i=0}^{l-1} \mathfrak{r}_{i}^{\sum\limits_{\mu \in \Z } N_{\mu l+i-j} (Y)}.
\end{gather*}
\end{Theorem}
\begin{proof}
Note that $\mathcal{Z}^{\rm{quot}}$ (resp.~$\mathcal{Z}^{\rm{core}}_{\mathbf{e}_j}$) is the contribution from $l$-quotients (resp.~$l$-cores) of~$\mathcal{P}_{\w}(\vb)$ corresponding to $T$-f\/ixed points on $\M_{\mathbf{e}_j} (\vb)$. By Corollary~\ref{cor:QforC}, we have
\begin{gather*}
\mathcal{Z}^{\mathrm{quot}} (\mathfrak{t},\mathfrak{q},\vec{\mathfrak{r}})= \sum_{n \geq 0} \mathfrak{q}^{nl} \prod_{i=0}^{l-1} \mathfrak{r}_{i}^{n}
			\sum_{Y \in \mathcal{P} (nl : \varnothing )} \mathfrak{t}^{ \# \{ s \in Y \,|\, h(s) \equiv 0 \, ({\rm{mod}} \, l ),\, l(s)>0 \} }.
\end{gather*}
Thus $\mathcal{Z}^{\rm{quot}} (\mathfrak{t},\mathfrak{q},\vec{\mathfrak{r}})$ is independent of $\mathbf{e}_j$. Then our claim follows from equation~(\ref{eq:v}) and the def\/inition of quotients and cores.
\end{proof}

\begin{Lemma}\label{lem:t2}
	\begin{gather*}
		\mathcal{Z}^{\rm{quot}} (\mathfrak{t},\mathfrak{q}_{\rm{reg}}) = \prod_{m=1}^{\infty} \frac{1}{\big(1-\mathfrak{q}_{\rm{reg}}^m \mathfrak{t}^{2m-2} \big)\big(1-\mathfrak{q}_{\rm{reg}}^m \mathfrak{t}^{2m} \big)^{l-1}} ,
	\end{gather*}
	where $\mathfrak{q}_{\rm{reg}} \defeq \mathfrak{q}^{l} \prod\limits_{i=0}^{l-1} \mathfrak{r}_{i}$.
\end{Lemma}
\begin{proof}
We have
\begin{gather*}
\mathcal{Z}^{\mathrm{quot}} (\mathfrak{t},\mathfrak{q}_{\rm{reg}} ) = \sum_{n \geq 0} \mathfrak{q}^{nl} \prod_{i=0}^{l-1} \mathfrak{r}_{i}^{n}
			\sum_{Y \in \mathcal{P} (nl : \varnothing )}
			\mathfrak{t}^{ \# \{ s \in Y \,|\, h(s) \equiv 0 \,({\rm{mod}} \, l ),\, l(s)>0 \} }\\
\hphantom{\mathcal{Z}^{\mathrm{quot}} (\mathfrak{t},\mathfrak{q}_{\rm{reg}} ) }{}
=\sum_{n \geq 0} \mathfrak{q}_{\rm{reg}}^{n}\left\{ \sum_{m\geq 0}\mathrm{b}_{m} (\M_{\mathbf{e}_{0}}(n\boldsymbol{\delta}) ) \mathfrak{t}^{m} \right\}.
\end{gather*}			
Thus $\mathcal{Z}^{\mathrm{quot}} (\mathfrak{t},\mathfrak{q}_{\rm{reg}} )$ is the generating function of the Poincar\'e polynomials of $\M_{\mathbf{e}_{0}}(n\boldsymbol{\delta})$. It is known~\cite{Kuznetsov, Wang} that $\M_{\mathbf{e}_0}(n\boldsymbol{\delta})$ is dif\/feomorphic to the Hilbert scheme of~$n$ points on the ALE space. By applying the G\"ottsche's formula~(\ref{eq:Gottsche}) to the ALE space, we obtain the result.
\end{proof}

\begin{Notation}\label{not:theta}
We use the following def\/inition of the Riemann theta function:
\begin{gather*}	
\boldsymbol\Theta (\mathbf{u} |\mathbf{F})\defeq \sum_{\mathbf{n} \in {\Z}^{l-1} }\exp \left[2\pi \sqrt{-1} \left({\mathbf{n}}^t \mathbf{u} +\frac{1}{2} \mathbf{n} \mathbf{F}^t \mathbf{n} \right)\right] ,
\end{gather*}
where $\mathbf{u}$ is a complex $(l-1)$-dimensional row vector, $^t (\, \, )$ means the transposition, and $\mathbf{F}$ is an $(l-1) \times (l-1)$ matrix.
\end{Notation}
\begin{Lemma}\label{lem:t3}We have
	\begin{gather*}
		\mathcal{Z}^{\rm{core}}_{\mathbf{e}_{j}} (\mathfrak{q},\vec{\mathfrak{r}}) = \boldsymbol\Theta (\mathbf{u}_j |\mathbf{F}_{\Gamma}),
	\end{gather*}
	where
	\begin{gather*}
		\mathbf{u}_j \defeq \frac{1}{2\pi \sqrt{-1}} (\log \mathfrak{q}\mathfrak{r}_{1+j} ,\ldots ,	\log \mathfrak{q}\mathfrak{r}_{{l-1}+j} )
	\end{gather*}
	and
	\begin{gather*}
		\mathbf{F}_{\Gamma} \defeq \frac{\log \mathfrak{q}_{\rm{reg}}}{2\pi \sqrt{-1}} \mathbf{C}_{\Gamma},
	\end{gather*}
	where $\mathbf{C}_{\Gamma}$ is the Cartan matrix of type $A_{l-1}$.
\end{Lemma}
\begin{proof}
	By Proposition \ref{prop:t1} and its proof, we have
	\begin{gather*}
\mathcal{Z}^{\rm{core}}_{\mathbf{e}_j} (\mathfrak{q},\vec{\mathfrak{r}}) =\sum_{ Y \in \mathcal{C}^{(l)} } \mathfrak{q}^{|Y|} \prod_{i=0}^{l-1} \mathfrak{r}_{i}^{\sum\limits_{\mu \in \Z } N_{\mu l+i-j} (Y)} \\
\hphantom{\mathcal{Z}^{\rm{core}}_{\mathbf{e}_j} (\mathfrak{q},\vec{\mathfrak{r}})}{}
	= \sum_{(k_{0}, \ldots, k_{l-1})\in \mathcal{K}^{(l)}}
				\mathfrak{q}^{\frac{1}{2} l \sum\limits_{i=0}^{l-1} k_i^2 +\sum\limits_{i=0}^{l-1}i k_i}
				\mathfrak{r}_{j}^{\frac{1}{2}\sum\limits_{i=0}^{l-1}k_{i}^2}
				\prod_{\mu =1}^{l-1}\mathfrak{r}_{j+\mu }^{\frac{1}{2}\sum\limits_{i=0}^{l-1}k_{i}^2 +\sum\limits_{i=\mu }^{l-1}k_{i}}\\
\hphantom{\mathcal{Z}^{\rm{core}}_{\mathbf{e}_j} (\mathfrak{q},\vec{\mathfrak{r}})}{}
			= \sum_{(k_{0}, \ldots, k_{l-1})\in \mathcal{K}^{(l)}}
				\mathfrak{q}_{\rm{reg}}^{\frac{1}{2} \sum\limits_{i=0}^{l-1} k_i^2 }
				\prod_{\mu =1}^{l-1} (\mathfrak{q}\mathfrak{r}_{j+\mu } )^{\sum\limits_{i=\mu }^{l-1} k_i} \\
\hphantom{\mathcal{Z}^{\rm{core}}_{\mathbf{e}_j} (\mathfrak{q},\vec{\mathfrak{r}})}{}
			= \sum_{( n_1, \ldots, n_{l-1} ) \in {\Z}^{l-1}}
				\prod_{\mu =1}^{l-1} \mathfrak{q}_{\rm{reg}}^{n_\mu (n_\mu -n_{\mu +1} )}
				(\mathfrak{q}\mathfrak{r}_{j+\mu } )^{n_\mu} ,
	\end{gather*}			
where $n_\mu \defeq \sum\limits_{i=\mu }^{l-1} k_i$. Comparing this equation with the theta function in Notation~\ref{not:theta}, we obtain the lemma.
\end{proof}

\begin{Theorem}\label{thm:t1}
We have
	\begin{gather*}
		\mathcal{Z}_{\mathbf{e}_{j}} (\mathfrak{t},\mathfrak{q},\vec{\mathfrak{r}})=\frac{\boldsymbol\Theta (\mathbf{u}_j |\mathbf{F}_{\Gamma})}{\prod\limits_{m=1}^{\infty} \big(1-\mathfrak{q}_{\rm{reg}}^m \mathfrak{t}^{2m-2} \big)\big(1-\mathfrak{q}_{\rm{reg}}^m \mathfrak{t}^{2m} \big)^{l-1}} .
	\end{gather*}
\end{Theorem}
\begin{proof}
This follows from Theorem \ref{thm:t2}, Lemmas~\ref{lem:t2} and~\ref{lem:t3}.
\end{proof}

\begin{Remark}\label{rem:nakajima} As we mentioned in the introduction, Theorem \ref{thm:t1} is obtained from the result of Nakajima \cite[Section~5.2]{Nakajima3} (see the comments below the remark). This was informed us by H.~Nakajima after the original version of this paper was submitted to e-print archives. He also pointed out that our argument has a close parallel to a geometric Frenkel--Kac construction of the Fock spaces of af\/f\/ine Lie algebras due to Grojnowski~\cite{Grojnowski} (see also \cite[Section~9.5]{Lecture}). The authors are grateful to him for these comments.
\end{Remark}

Let us make some comments on representation theoretical aspects of the result in Theorem~\ref{thm:t1}. By Nakajima's pioneering work, it is known that quiver varieties are deeply related to representations of (quantum) af\/f\/ine algebras $\mathfrak{g}$ (of type $A_{l-1}^{(1)}$ in our case). It is shown in~\cite{Nakajima4,Nakajima6} that their ($K$-)homology groups have structures of representations of af\/f\/ine algebras. Roughly speaking, $\oplus_{\vb} H_{*} (\M_{\w}(\vb) )$ is a direct sum of certain highest weight representations of $\mathfrak{g}$ determined by~$\w$. When $\w$ is the class of trivial $1$-dimensional $\Gamma$-module, this is the so-called \textit{Fock space} representation of~$\mathfrak{g}$. This is explained in \cite[Section~5.2]{Nakajima3}. From this view point, the generating function of Euler characteristics is the character of the representation and the
generating function of Poincar\'e polynomials is the so-called \textit{q,t-character}~\cite{Nakajima5}, which plays a~fundamental role in representation theory.

\subsubsection{Euler characteristics for higher rank case}
Let $\dim W=r$. We consider the following generating function of the Euler characteristics ${\rm e} ( \M_{\w} (\vb))$ of $\M_{\w} (\vb)$:
\begin{gather*}
\mathcal{E}_{\w} (\mathfrak{q}) \defeq \sum_{\vb} {\rm e} ( \M_{\w} (\vb)) \mathfrak{q}^{|\vb |}\prod_{i=0}^{l-1} \mathfrak{r}_i ^{v_i }.
\end{gather*}
\begin{Corollary}
	\begin{gather*}
\mathcal{E}_{\w} (\mathfrak{q})=\mathfrak{q}_{\rm{reg}}^{\frac{rl}{24}}	\frac{\prod\limits_{j=0}^{l-1} \boldsymbol\Theta (\mathbf{u}_j |\mathbf{F}_{\Gamma})^{w_j}}{\eta (\mathfrak{q}_{\rm{reg}} )^{rl}},
	\end{gather*}
	where $\eta (q) \defeq q^{\frac{1}{24}} \prod\limits_{n=1}^{\infty} (1- q^n)$ is the Dedekind eta function.
\end{Corollary}
\begin{proof}
We have
\begin{gather*}
	\mathcal{E}_{\w} (\mathfrak{q})=\prod_{j=0}^{l-1} (\mathcal{Z}_{\mathbf{e}_{j}} (\mathfrak{t}=1,\mathfrak{q},\vec{\mathfrak{r}}) )^{w_j} .
\end{gather*}
Then by Theorem~\ref{thm:t1}, we have the claim.
\end{proof}

\section{Relations to instanton counting} \label{section5}
In this section, we discuss connections with \textit{instanton counting} in physics. In Section~\ref{sec:IC}, we give a brief explanation of instanton counting, which is one of our motivations of the f\/irst part of this paper. In Section~\ref{section5.2}, we explain a part of results of Maeda et al.~\cite{Nakatsu3, Nakatsu1,Nakatsu2}, where quotients and cores for Young diagrams appeared.

\subsection{Instanton counting}\label{sec:IC}
Mathematically speaking, instanton counting means computations of global invariants on the instanton moduli space, such as the Euler characteristic and the volume of the moduli space. What kind of invariants one want to compute depends on what kind of physical theory one considers. As we mentioned in introduction, $\M_{\w}(\vb)$ is a resolution of the moduli space of instantons on $\C^2/\Gamma$, therefore our results in Section~\ref{sec:enu} can be regarded as an instanton counting on $\C^2/\Gamma$.

\subsubsection{}
In four-dimensional ($4$D) topologically twisted $\mathcal{N}=4$ supersymmetric (SUSY) Yang--Mills (YM) theory, the partition function is given by the generating function of \textit{Euler characteristics} of instanton moduli spaces, as we considered in Section~\ref{sub:euler}. Such a theory is studied extensively by Vafa and Witten \cite{Vafa-Witten}. See also \cite{Jinzenji,Sasaki}, which considered issues similar to ours.

\subsubsection{}\label{sub:n=2}
Nekrasov \cite{Nek1} introduced
the following generating function of \textit{equivariant volumes} of the instanton moduli spaces, in the sense of Section~\ref{sec:integral}, in studies of $4$D $\mathcal{N}=2$ SUSY YM theory:
\begin{gather}
Z^{\mathrm{inst}}({\epsilon}_1 ,{\epsilon}_2 , \vec{a}; \Lambda)\defeq \sum_{n=0}^{\infty}\Lambda^{2nr} \int_{\M_{r}(n)} 1 .\label{eq:nekpar}
\end{gather}
Such a generating function is called (the instanton part of) Nekrasov's partition function (see also~\cite{BFMT, FP}). Nekrasov conjectured that $F^{\mathrm{inst}}({\epsilon}_1 ,{\epsilon}_2 , \vec{a}; \Lambda) \defeq \epsilon_1\epsilon_2 \log Z^{\mathrm{inst}}({\epsilon}_1 ,{\epsilon}_2 , \vec{a}; \Lambda)$ is regular at $\epsilon_1=\epsilon_2=0$, and $F^{\mathrm{inst}}(0, 0, \vec{a}; \Lambda)$ coincides with (the instanton part of) the celebrated Seiberg--Witten prepotential. This conjecture is proved independently by Nakajima and Yoshioka \cite{NY2, NY}, Nekrasov and Okounkov~\cite{NekOk}, and
Braverman and Etingof \cite{Braverman1, Braverman2}. We do not give a detail about Seiberg--Witten theory here. We refer the mathematical oriented readers \cite[Section~2]{NY2} for the Seiberg--Witten prepotential.

\subsubsection{}\label{sub:5d}
There is a $5$D version $Z^{\mathrm{inst}}_{5\mathrm{D}}({\epsilon}_1 ,{\epsilon}_2 , \vec{a}; \Lambda)$ of the Nekrasov's partition function~(\ref{eq:nekpar}), which was also introduced by Nekrasov \cite{Nek1}. Physically speaking, that is $\mathcal{N}=1$ SUSY YM theory on $\R^4 \times S^1$ (in a certain supergravity background). Mathematically speaking, this means that we consider equivariant integrals in (\ref{eq:nekpar}) in the sense of equivariant $K$-theory \cite{NY3}.

\subsection{Cores and perturbative gauge theory}\label{section5.2}
In this section, we consider $5$D version of the Nekrasov's partition function and discuss a~combinatorial aspect of it. We explain works of Maeda et al.~\cite{Nakatsu3, Nakatsu1,Nakatsu2}, where quotients and cores for Young diagrams appeared in a study of $5$D $\mathcal{N}=1$ SUSY YM theory explained in Section~\ref{sub:5d}.

\subsubsection{}
So far, we consider only the instanton part of the partition function. There is the so-called perturbative part $Z^{\mathrm{pert}}_{5\mathrm{D}}({\epsilon}_1 ,{\epsilon}_2 , \vec{a}; \Lambda)$ of the partition function, which corresponds to the perturbative part of the Seiberg--Witten prepotential. Although we do not give a def\/inition here, it is def\/ined by an explicit formula (see, for example,~\cite[Appendix~E]{NY2}). Then the full partition function is def\/ine by
\begin{gather*}
Z_{5\mathrm{D}}({\epsilon}_1 ,{\epsilon}_2 , \vec{a}; \Lambda) = Z^{\mathrm{pert}}_{5\mathrm{D}}({\epsilon}_1 ,{\epsilon}_2 , \vec{a}; \Lambda)\, Z^{\mathrm{inst}}_{5\mathrm{D}}({\epsilon}_1 ,{\epsilon}_2 , \vec{a}; \Lambda).%\label{eq:factorization}
\end{gather*}

\subsubsection{Results of Maeda et al.}
Recall that $T$-f\/ixed points on $\M_r (n)$ are parametrized by $\mathcal{P}_r (n)$ and an equivariant integral (in the sense of $K$-theory) over $\M_r (n)$ can be expressed as a sum over $\mathcal{P}_r (n)$ by localization theorem, as we discussed in Section~\ref{sec:integral}. The authors of  loc.\ cit.\ found the following identity:
\begin{gather}
Z^{r=1}_{5\mathrm{D}}(\hbar, -\hbar ; \Lambda) =\sum_{\vec{k}\in \mathcal{K}^{(N)}} Z^{r=N}_{5\mathrm{D}}\big(\hbar, -\hbar,
\vec{\tilde{k}}; \Lambda \big),\label{eq:nakatsu}
\end{gather}
where the LHS of (\ref{eq:nakatsu}) is the partition function for $r=1$ under the condition\footnote{Under this condition, Nekrasov's partition function can be expressed in terms of special values of Schur functions.}
\begin{gather}
\epsilon_1=-\epsilon_2=\hbar\label{eq:schur}
\end{gather}
and summands in the RHS of
(\ref{eq:nakatsu}) are the partition function for $r=N(\geq 2)$ under the condition~(\ref{eq:schur}) and specialization $\vec{a} = \vec{\tilde{k}}$, where $\vec{\tilde{k}} =(\tilde{k}_{1}, \ldots , \tilde{k}_{N})$ is given by
\begin{gather*}
\tilde{k}_{m}=\hbar \left\{ k_{m}+\frac{1}{N}\left(m-\frac{N+1}{2}\right) \right\}
\end{gather*}
for $\vec{k}=(k_1, \ldots , k_N) \in \mathcal{K}^{(N)}$. In  loc.\ cit., this identity was proved by the following way. First, $Z_{5\mathrm{D}}(\hbar, -\hbar, \vec{\tilde{k}}; \Lambda)$ can be identif\/ied with the partition function of a statistical model of plane partitions (= $3$D Young diagrams) whose $r$-core of the `main diagonal' is $Y^{(r)} (\vec{k})$, which can be described by free fermions. Then, under this identif\/ication, the perturbative part $Z^{\mathrm{pert}}_{5\mathrm{D}}(\hbar, -\hbar, \vec{\tilde{k}}; \Lambda)$ is identif\/ied with the contribution of ground states, which correspond to $N$-cores $Y^{(N)}(\vec{k}) \in \mathcal{C}^{(N)}$ determined by $\vec{k} \in \mathcal{K}^{(N)}$. And the instanton part $Z^{\mathrm{inst}}_{5\mathrm{D}}(\hbar, -\hbar, \vec{\tilde{k}}; \Lambda)$ is identif\/ied with the contribution of excitations, which correspond to $N$-quotients. Notice that, by Proposition~\ref{thm:s3.6}, we recover the set of all Young diagrams $\mathcal{P}$ if we vary all $\vec{k}=(k_1, \ldots, k_N) \in \mathcal{K}^{(N)}$ and $\underline{Y}\in\mathcal{P}_{N}$. Taking these into account, we can see that the identity~(\ref{eq:nakatsu}) holds, after some Schur function calculus.

It is worth mentioning that we have a similar factorization of the generating function $\mathcal{Z}_{\mathbf{e}_{j}}\! (\mathfrak{t}{,} \mathfrak{q}{,}\vec{\mathfrak{r}})$ of Poincal\'e polynomials in rank $1$ case into `core-part' $\mathcal{Z}^{\rm{core}}$ and `quotient-part' $\mathcal{Z}^{\rm{quot}}_{\mathbf{e}_j}$ in Theo\-rem~\ref{thm:t2}. Thus combinatorial structures in our case and the case studied by Maeda et al.\ are essentially the same.
\begin{Remark} Under the condition~(\ref{eq:schur}), it is known that Nekrasov's partition function has connections with topological strings \cite{Eguchi-Kanno, Iqbal, Tachikawa, Zhou} and $2$D YM theory \cite{Matsuura-Ohta}. See Remark~\ref{rem:mac} for a~comment on cases with $\epsilon_1 + \epsilon_2 \neq 0$.
\end{Remark}

\appendix
\section{Hilbert schemes of points on surfaces} \label{appendixA}
This appendix has two purposes. The one is to summarize some facts about the Hilbert schemes of points\footnote{A basic reference on this subject is \cite{Lecture}.} (Section~\ref{sub:app}). The other is to give another method to compute equivariant volumes of quiver varieties in rank $1$ case (Sections~\ref{sub:torus} and~\ref{sec:example}). We get a closed formula for the generating function of such quantities. By comparing with the result in Section~\ref{sec:integral}, we f\/ind an interesting combinatorial identity (Proposition~\ref{prop:s4.6}).

\subsection{Hilbert schemes of points}\label{sub:app}

\subsubsection{}
Let $X$ be a smooth quasi-projective surface. The Hilbert scheme $\Hilbn{X}$ of $n$ points on $X$ is, by def\/inition, a smooth algebraic variety of dimension $2n$, which parametrizes the set of 0-dimensional subschemes in~$X$ with colength~$n$. Topology of the Hilbert scheme $\Hilbn{X}$ is well-known. The generating function of the Poincar\'e polynomials is given by the following G\"ottsche's formula~\cite{gottsche}:
	\begin{gather}
		\sum_{n\geq 0, i\geq 0} b_i \big(\Hilbn{X} \big)\mathfrak{t}^i \mathfrak{q}^n
			=\prod_{m \geq 1, i\geq 0} \big( 1-(-1)^i \mathfrak{t}^{2m-2+i} \mathfrak{q}^m \big)^{(-1)^{i+1} b_i (X)}.	\label{eq:Gottsche}
	\end{gather}

\subsubsection{}
\label{sec:action}
There exists a proper surjective morphism $\pi \colon \Hilbn{X} \to S^{n}X$, where $S^{n}X$ is the $n$-th symmetric product of $X$. This morphism, so called the \textit{Hilbert--Chow morphism}, gives a resolution of singularities of $S^{n}X$. Note that if $X$ admits a torus action, then there is an induced torus action on $\Hilbn{X}$ and the Hilbert--Chow morphism $\pi$ is equivariant.

\subsection{Torus actions on Hilbert schemes}\label{sub:torus}
We study torus actions on the Hilbert schemes of points on toric surfaces. All the results in this section are due to Ellingsrud--Str{\o}mme \cite{ES} (see also \cite{Li1,Li2,Qin}).

\subsubsection{Toric surfaces}
Let $N$ be a 2-dimensional lattice and $\Sigma$ be a fan in $N_{\R}$. From these datum, we can construct an algebraic surface $X = X(N, \Sigma)$, equipped with an action of the 2-dimensional algebraic torus $T=(\C^{*})^2$ with a dense open orbit. These surfaces are called toric surfaces, and it is known that they are normal and quasi-projective (see, for example, \cite{Fulton, Oda} for detail).

\subsubsection{Torus actions on toric surfaces}
Let $X$ be a smooth toric surface associated to $(N, \Sigma)$. Then the $2$-dimensional torus $T=N\otimes_{\Z}\C^{*}$ acts on~$X$. The $T$-f\/ixed points correspond to vertices of $2$-dimensional cones in~$\Sigma$. In particular, the number of $T$-f\/ixed points is equal to the Euler number ${\rm e}(X)$ of~$X$. Let $X^{T} = \{p_1, \ldots, p_{{\rm e}(X)} \}$. Then, by construction, there is an af\/f\/ine chart $U_i =\operatorname{Spec} \C[x_i, y_i]$ around each $p_i$ which is $T$-invariant. We call $(x_i, y_i)$ a toric coordinate around $p_i$. We denote $w_{x_i}(t_1, t_2)$ (resp.~$w_{y_i}(t_1, t_2)$) the weight of the $T$-action on $x_i$ (resp.~$y_i$).

\subsubsection{Torus actions on Hilbert schemes}
As we noticed in Section~\ref{sec:action}, $T$ acts also on $\Hilbn{X}$. First, we identify $T$-f\/ixed points on $\Hilbn{X}$.
\begin{Lemma}
There is a one-to-one correspondence between $T$-fixed point set on $\Hilbn{X}$ and $\mathcal{P}_{{\rm e}(X)}(n)$.
\end{Lemma}
\begin{proof}
A 0-dimensional subscheme $Z$ of length $n$ on $X$ is $T$-f\/ixed if and only if $Z$ is supported on $T$-f\/ixed points $\left\{p_i\right\}$ on~$X$. Each connected component $Z_i$ of $Z$ with $\Supp(Z_i) = \{p_i\}$ def\/ines a $T$-f\/ixed point in $(U_i)^{[n_i]} \cong \big(\C^2\big)^{[n_i]}$ for some $n_i (\geq 0)$ with $\sum\limits_{i=1}^{{\rm e}(X)} n_i =n$. It is well-known that there is a one-to-one correspondence between $T$-f\/ixed points in $\Hilbn{\big(\C^2\big)}$ and partitions of~$n$ (see \cite[Chapter~5]{Lecture}). This gives the assertion.
\end{proof}

\begin{Lemma} \label{eq:char}
The weight decomposition of the cotangent space of $\Hilbn{X}$ at a $T$-fixed point corresponding to $\underline{Y} =( Y_{1}, \ldots, Y_{{\rm e}(X)}) \in \mathcal{P}_{{\rm e}(X)}(n)$ is given by
\begin{gather*}
 \sum_{i=1}^{{\rm e}(X)} \sum_{s \in Y_{i}} \big\{ w_{x_i}(t_1, t_2)^{l(s)+1} w_{y_i}(t_1, t_2)^{-a(s)} + w_{x_{i}}(t_1, t_2)^{-l(s)} w_{y_i}(t_1, t_2)^{a(s)+1} \big\}.
\end{gather*}
\end{Lemma}

 \begin{proof}
 This is obvious from the classical result of Ellingsrud and Str{\o}mme \cite{ES} on the weight decomposition of the tangent space of $\Hilbn{\big(\C^2\big)}$ at
 $T$-f\/ixed point (see \cite[Proposition~5.8]{Lecture}).
\end{proof}

For a Young diagram $Y$, we set
\begin{gather*}
n_Y (\eta_1 , \eta_2) \defeq \prod_{s\in Y}\{ ( l(s)+1) \eta_1 -a(s) \eta_2\}\{-l(s) \eta_1 +( a(s)+1) \eta_2 \}.
\end{gather*}
Here $\eta_1$, and $\eta_2$ are parameters. Then we have the following
\begin{Proposition}
$T$-equivariant integral on $X^{[n]}$ is given by
\begin{gather*}
\int_{\Hilbn{X}} 1 = \sum_{\underline{Y} \in \mathcal{P}_{{\rm e}(X)}(n)} \prod_{i=i}^{{\rm e}(X)} \frac{1}{n_{Y_i}
\big( \log w_{x_i}\big(e^{\epsilon_1}, e^{\epsilon_2}\big), \log w_{y_i}(e^{\epsilon_1}, e^{\epsilon_2})\big)}.
\end{gather*}
Here we write $t_i = e^{\epsilon_i}$, where $\epsilon_i \in \operatorname{Lie}(T)$.
\end{Proposition}

\subsubsection{}\label{sub:gen}
We consider the following generating function (\ref{eq:nekpar}):
\begin{gather*}
Z_{X}(\epsilon_1 , \epsilon_2 ; \Lambda) \defeq \sum_{n=0}^{\infty} \Lambda^{2n} \int_{\Hilbn{X}} 1 .
\end{gather*}

Let us consider the case $X=\C^2$ equipped with a $T$-action def\/ined by
\begin{gather}\label{eq:action}
(x,y) \mapsto (t_{1}x, t_{2}y).
\end{gather}
It is shown in \cite[Section~4]{NY} that
\begin{gather*}
Z_{\C^2}(\epsilon_1 , \epsilon_2 ; \Lambda)= \exp \left( \frac{\Lambda^2}{\epsilon_1 \epsilon_2} \right).
\end{gather*}
Then we have
\begin{gather}%\label{eq:product}
Z_{X} (\epsilon_1 , \epsilon_2 ; \Lambda) =\sum_{n=0}^{\infty} \sum_{\underline{Y} \in \mathcal{P}_{{\rm e}(X)}(n)}
\frac{\Lambda^{2|\underline{Y}|}}{\prod\limits_{i=i}^{{\rm e}(X)} n_{Y_i}
\big( \log w_{x_i}\big(e^{\epsilon_1}, e^{\epsilon_2}\big), \log w_{y_i}(e^{\epsilon_1}, e^{\epsilon_2})\big)}\nonumber\\
\hphantom{Z_{X} (\epsilon_1 , \epsilon_2 ; \Lambda)}{} = \prod_{i=1}^{{\rm e}(X)} Z_{\C^2}\big( \log w_{x_i}\big(e^{\epsilon_1}, e^{\epsilon_2}\big), \log w_{y_i}(e^{\epsilon_1}, e^{\epsilon_2}) ; \Lambda\big)\nonumber\\
\hphantom{Z_{X} (\epsilon_1 , \epsilon_2 ; \Lambda)}{}= \exp \left\{\Lambda^2\sum_{i=1}^{{\rm e}(X)} \frac{1}
{\big( \log w_{x_i}(e^{\epsilon_1}, e^{\epsilon_2})\big) \big( \log w_{y_i}(e^{\epsilon_1},
e^{\epsilon_2})\big)}\right\}. \label{eq:formula}
\end{gather}

\subsection{An example}\label{sec:example}
We study a particular case when $X$ is the ALE space of type $A_{l-1}$, i.e., $X=\M_{{\rm e}_0}(\boldsymbol\delta)$ and compare the result in Section~\ref{sub:gen} and the result in Section~\ref{sec:integral}.

\subsubsection{}
Recall that the simple singularity of type $A_{l-1}$ is the quotient space $\C^2/\Gamma$, where the action of $\Gamma$ on $\C^2$ is given by $(x, y) \mapsto (\zeta x, \zeta ^{-1} y)$. Here $\zeta$ is a primitive $l$-th root of unity. Note that the above action of $\Gamma$ on $\C^2$ commutes with the $T$-action~(\ref{eq:action}) on~$\C^2$. It follows that $\C^2/\Gamma$ is a toric singularity. It is well-known \cite{Fulton, Oda} that it has the toric minimal resolution. In~\cite{IN1}, it is shown that $X$ is isomorphic to the toric minimal resolution, so that~$X$ has af\/f\/ine charts $U_i =\mathrm{Spec}\,\C[x_i, y_i]$ $(i=1, \ldots ,l)$ def\/ined by
\begin{gather}
x_i = \frac{a^i}{b^{l-i}}, \qquad y_i = \frac{b^{l+1-i}}{a^{i-1}},\label{eq:coord}
\end{gather}
where $(a, b)$ is a coordinate of $\C^2$ on which $\Gamma$ acts by $(a,b) \mapsto (\zeta a, \zeta ^{-1} b)$.
\begin{Remark}
It is easy to see that there is a bijection between the set of these $l$ af\/f\/ine charts and the set $\mathcal{P}_{\mathbf{e}_0}(\boldsymbol\delta)$. This correspondence is pointed out in \cite[Corollary~3.10]{Yukarin}, where $\mathcal{P}_{\mathbf{e}_0}(\boldsymbol\delta)$ is identif\/ied with the set of so-called $\Gamma$-\textit{clusters}.
\end{Remark}

From equation (\ref{eq:coord}), we have
\begin{gather*}
\log w_{x_{i}} (\epsilon_1, \epsilon_2) = i\epsilon_1 -(l-i)\epsilon_2 ,\\
\log w_{y_{i}} (\epsilon_1, \epsilon_2) = -(i-1)\epsilon_1 +(l+1-i)\epsilon_2 .
\end{gather*}
By substituting these into (\ref{eq:formula}), we get
\begin{gather*}
Z_{X} (\epsilon_1 , \epsilon_2 ; \Lambda) = \exp \left( \frac{\Lambda^2}{l\epsilon_1 \epsilon_2}\right) .
\end{gather*}
The coef\/f\/icient of $\Lambda^{2n}$ is given by
\begin{gather*}
\frac{1}{n!l^n (\epsilon_{1} \epsilon_{2})^{n}},%\label{eq:B9}
\end{gather*}
which is nothing but the $T$-equivariant volume of the orbifold $S^n\big(\C^2/\Gamma\big)$. This is the desired result, since natural morphism $\Hilbn{X} \to S^n\big(\C^2/\Gamma\big)$, induced by $X \to \C^2/\Gamma$, is $T$-equivariant.

\subsubsection{}
We have another resolution of singularities of $S^n\big(\C^2/\Gamma\big)$, which is given by quiver varieties \cite[Section~7.2.3]{Haiman2} (see also~\cite{Wang}). Note that by Propositions \ref{thm:s3.6} and~\ref{prop:s3.8}, total weight~$|\underline{Y}_{(l)}^*|$ of $l$-quotient of $Y \in \mathcal{P}_{\mathbf{e}_j}(\vb)$ is independent of~$Y$. Let $n = |\underline{Y}_{(l)}^*|$ for $Y \in \mathcal{P}_{\mathbf{e}_j}(\vb)$. Then by using Nakajima's \textit{reflection functor} \cite{Nakajima7}, it can be shown that $\M_{\mathbf{e}_j } (\vb)$ is a $T$-equivariant resolution of singularities of $S^n\big(\C^2/\Gamma\big)$. Thus equivariant integral $\int_{\M_{\mathbf{e}_j} (\vb)} 1$ also gives the equivariant volume of~$S^n\big(\C^2/\Gamma\big)$. Therefore we have the following identity.
\begin{Proposition}\label{prop:s4.6}
	\begin{gather*}
	 \sum_{Y \in \mathcal{P}_{\mathbf{e}_j} (\vb )}	\prod_{s \in Y \atop h (s) \equiv 0 \, {\rm{mod}} \, l}	\frac{1}{( -l (s) {\epsilon}_1 +(a (s)+1){\epsilon}_2 )( (l (s)+1) {\epsilon}_1 -a (s){\epsilon}_2 )} = \frac{1}{n! l^n ({\epsilon}_1 {\epsilon}_2 )^n} ,
	\end{gather*}
	where $n = |\underline{Y}_{(l)}^*|$ for $Y \in \mathcal{P}_{\mathbf{e}_j}(\vb)$.
\end{Proposition}
Since the above identity is purely combinatorial, it should be possible to give a combinatorial proof. In the case of ${\epsilon}_1 = -{\epsilon}_2 =\hbar$, we have the following combinatorial proof.
\begin{gather*}
\sum_{Y \in \mathcal{P}_{\mathbf{e}_j} (\vb ) } \prod_{s \in \underline{Y}_{(l)}^* } \frac{1}{-( lh(s) \hbar)^2 }
		 = \frac{1}{l^{2n} (-{\hbar}^2 )^n}	\sum_{(n_0 ,\ldots ,n_{l-1} )\in {\Z}_{\geq 0}^{l} \atop n_0 + \cdots + n_{l-1} =n }
			\prod_{i=0}^{l-1} \sum_{Y \in P(n_i )} \prod_{s \in Y} \frac{1}{h(s)^2} \\
\hphantom{\sum_{Y \in \mathcal{P}_{\mathbf{e}_j} (\vb ) } \prod_{s \in \underline{Y}_{(l)}^* } \frac{1}{-( lh(s) \hbar)^2 }}{} = \frac{1}{l^{2n} (-{\hbar}^2 )^n}	\sum_{(n_0 ,\ldots , n_{l-1} )\in {\Z}_{\geq 0}^{l} \atop n_0 + \cdots + n_{l-1} =n }
			\prod_{i=0}^{l-1} \frac{1}{n_i !} = \frac{1}{n! l^n (-{\hbar}^2 )^n}.
\end{gather*}
However, the authors do not have a combinatorial proof when $\epsilon_1 + \epsilon_2 \neq 0$.
\begin{Remark}\label{rem:mac}
It is remarkable that in \cite{Haiman1}, geometry of the Hilbert scheme is applied to problems of enumerative combinatorics which is related to Macdonald polynomials. Note also that in a physical literature \cite{AwaKan}, relations between Macdonald polynomials and instanton counting in $5$D SUSY YM theory (see Section~\ref{sub:n=2}) is pointed out. It seems interesting to study further in this direction.
\end{Remark}

\subsection*{Acknowledgements}
The authors would like to thank
H.~Awata, H.~Miyachi, W.~Nakai, H.~Nakajima, T.~Nakatsu, M.~Namba, Y.~Nohara, Y.~Hashimoto, Y.~Ito, T.~Sasaki, Y.~Tachikawa, K.~Takasaki, and K.~Ueda for valuable discussions and comments. The authors express their deep gratitudes to M.~Hamanaka, S.~Moriyama, and A.~Tsuchiya for their advices and warm encouragements, and especially to H.~Kanno for suggesting a problem and reading the manuscript carefully. This work was started while the authors enjoyed the hospitality of the Fields Institute at University of Toronto on the fall of 2004. The authors are grateful to K. Hori for invitation. Throughout this work, the authors' research was supported in part by COE program in mathematics at Nagoya University.

{\it Added in 2017.} The authors thank the referees for useful comments. During the revision in 2017, S.M.~is supported in part by Grant for Basic Science Research Projects from the Sumitomo Foundation and JSPS KAKENHI Grand number JP17K05228.

%\cite{Date,FrenkelSavage,Fujii-Minabe,Kac,Miyachi,Nakajima7}

\pdfbookmark[1]{References}{ref}
\LastPageEnding

\end{document}